\documentclass{ert-l}
\usepackage{amssymb,enumerate}
\newcommand{\FF}{{\mathbb{F}}}

\newcommand{\ZZ}{{\mathbb{Z}}}

\newcommand{\fS}{{\mathfrak{S}}}

\newcommand{\bG}{{\mathbf{G}}}
\newcommand{\bH}{{\mathbf{H}}}
\newcommand{\bL}{{\mathbf{L}}}
\newcommand{\bT}{{\mathbf{T}}}

\newcommand{\cH}{{\mathcal{H}}}
\newcommand{\cS}{{\mathcal{S}}}

\newcommand{\Irr}{{\operatorname{Irr}}}
\newcommand{\df}{{\operatorname{def}}}
\newcommand{\diag}{{\operatorname{diag}}}
\newcommand{\rk}{{\operatorname{rk}}}

\newcommand{\PSL}{{\operatorname{L}}}
\newcommand{\SL}{{\operatorname{SL}}}
\newcommand{\SU}{{\operatorname{SU}}}
\newcommand{\PSU}{{\operatorname{U}}}
\newcommand{\Sp}{{\operatorname{Sp}}}
\newcommand{\Spin}{{\operatorname{Spin}}}
\newcommand{\OO}{{\operatorname{O}}}
\newcommand{\SO}{{\operatorname{SO}}}
\newcommand\RLG{{R_\bL^\bG}}
\newcommand\sRLG{{\,{}^*R_\bL^\bG}}

\newcommand{\ul}[1]{\underline{#1}}
\newcommand{\sy}[2]{\ensuremath{\left(\begin{array}{c}{#1}\\{#2}\end{array}\right)}}

\newcommand{\tw}[1]{{}^#1\!}

\let\eps=\epsilon
\let\la=\lambda
\let\sm=\setminus
\let\sys=\binom

\newtheorem{thm}{Theorem}[section]
\newtheorem{lem}[thm]{Lemma}

\newtheorem{cor}[thm]{Corollary}
\newtheorem{prop}[thm]{Proposition}

\newtheorem{thmA}{Theorem}

\theoremstyle{definition}

\theoremstyle{remark}
\newtheorem{rem}[thm]{Remark}

\begin{document}

\title[Murnaghan--Nakayama rule for unipotent characters]
      {A Murnaghan--Nakayama rule for values\\ of unipotent characters in
       classical groups}

\date{\today}

\author{Frank L\"ubeck}
\address{Lehrstuhl D f\"ur Mathematik, RWTH Aachen, Pontdriesch 14/16,
  52062 Aachen, Germany.}
\email{Frank.Luebeck@math.rwth-aachen.de}

\author{Gunter Malle }
\address{FB Mathematik, TU Kaiserslautern, Postfach 3049,
         67653 Kaisers\-lautern, Germany.}
\email{malle@mathematik.uni-kl.de}

\thanks{The second author gratefully acknowledges financial support by ERC
  Advanced Grant 291512.}



\begin{abstract}
We derive a Murnaghan--Nakayama type formula for the values of unipotent
characters of finite classical groups on regular semisimple elements. This
relies on Asai's explicit decomposition of Lusztig restriction. We use our
formula to show that most complex irreducible characters vanish on some
$\ell$-singular element for certain primes $\ell$.

As an application we classify the simple endotrivial modules of the finite
quasi-simple classical groups. As a further application we show that for
finite simple classical groups and primes $\ell\ge3$ the first Cartan invariant
in the principal $\ell$-block is larger than~2 unless Sylow $\ell$-subgroups
are cyclic.
\end{abstract}

\maketitle

\section{Introduction}   \label{sec:intro}
The classical Murnaghan--Nakayama rule provides an efficient recursive method
to compute the values of irreducible characters of symmetric groups. This
method can be adapted to the finite general linear and unitary groups, where
it then allows to compute values of unipotent characters not on arbitrary
but just on regular semisimple elements, see \cite[Prop.~3.3]{LM15}. This
adaptation uses the fact that the unipotent characters of general linear groups
coincide with Lusztig's so-called almost characters and then applies the
Murnaghan--Nakayama rule for the symmetric group. \par
In the present paper, we derive a Murnaghan--Nakayama rule for the values of
unipotent characters of the finite classical groups on regular semisimple
elements, see Theorem~\ref{thm:MN}. Our approach here is not via a
corresponding formula for the underlying Weyl group (of type $B_n$ or $D_n$),
since for classical groups Lusztig's Fourier transform is non-trivial and thus
the relation between almost characters and unipotent characters becomes quite
tricky. Instead we combine Asai's result expressing the decomposition of
Lusztig induction in terms of hooks and cohooks of symbols and which already
has a Murnaghan--Nakayama like form, with Lusztig's character formula.
\medskip

As a first application we derive a vanishing result for irreducible characters
of quasi-simple groups on $\ell$-singular elements:

\begin{thmA}   \label{thmA:zeroes}
 Let $\ell>2$ be a prime and $G$ a finite quasi-simple group of $\ell$-rank at
 least~3. Then for any non-trivial character $\chi\in\Irr(G)$ there exists an
 $\ell$-singular element $g\in G$ with $\chi(g)=0$, unless either $G$ is of
 Lie type in characteristic~$\ell$, or $\ell=5$ and one of:
 \begin{enumerate}
  \item[\rm(1)] $G=\PSL_5(q)$ with $5||(q-1)$ and $\chi$ is unipotent
   of degree $\chi(1)=q^2\Phi_5$;
  \item[\rm(2)] $G=\PSU_5(q)$ with $5||(q+1)$ and $\chi$ is unipotent
   of degree $\chi(1)=q^2\Phi_{10}$;
  \item[\rm(3)] $G=Ly$ and $\chi(1)\in\{48174,11834746\}$; or
  \item[\rm(4)] $G=E_8(q)$ with $q$ odd, $d_\ell(q) = 4$ and $\chi$ one 
   character in the Lusztig-series of type $D_8$.
 \end{enumerate}
\end{thmA}

The second main ingredient in the proof, apart from our Murnaghan--Nakayama
formula, is a result asserting the existence of $\ell$-singular regular
semisimple elements in classical groups of Lie type (see
Lemma~\ref{lem:regelts}) which may be of independent interest.
\medskip

We then apply this vanishing result in order to classify simple endotrivial
modules over fields $k$ of positive characteristic $\ell$ for finite classical
groups.

\begin{thmA}   \label{thmA:endoclass}
 Let $G$ be a finite quasi-simple group of classical Lie type $B_n,C_n,D_n$ or
 $\tw2D_n$ with non-cyclic Sylow $\ell$-subgroup. Then there exists a
 non-trivial simple endotrivial $kG$-module $V$ if and only if $G=\Sp_8(2)$,
 $k$ has characteristic~$\ell=5$, and $\dim(V)=51$.
\end{thmA}

Our approach relies on the fact, proven in \cite[Thm.~1.3]{LMS13} that any
endotrivial module is liftable to a characteristic~0 representation, which can
then be studied by ordinary character theory. In particular, its character
cannot vanish on any $\ell$-singular element, and this latter condition can be
checked with the previous Murnaghan--Nakayama formula.
\medskip

As a further application we give a (partial) answer to a question of
Koshitani, K\"uls\-hammer and Sambale \cite{KKS}, which is then used in
\cite{LM15} to settle this question completely:

\begin{thmA}   \label{thmA:PIMs}
 Let $G$ be a finite simple group of classical Lie type $B_n,C_n,D_n$ or
 $\tw2D_n$. Let $\ell>2$ be a prime for which Sylow $\ell$-subgroups of $G$
 are non-cyclic. Then the $\ell$-modular projective cover of the trivial
 character of $G$ has at least three ordinary constituents. \par
 In particular the first Cartan invariant of $G$ satisfies $c_{11}\ge3$.
\end{thmA}

Our paper is organised as follows. In Section~\ref{sec:maxtor} we prove the
existence of regular semisimple elements in suitable maximal tori of
groups of classical type and determine possible overgroups of collections of
maximal tori. Section~\ref{sec:MNrule} contains the proof of the
Murnaghan--Nakayama rule, see Theorem~\ref{thm:MN}. We apply this in
Section~\ref{sec:vanish} to prove a vanishing result on $\ell$-singular
elements, see Theorems~\ref{thm:non-unip} and~\ref{thm:unip}. In
Section~\ref{sec:endotriv} we show Theorem~\ref{thmA:endoclass}
classifying simple endotrivial modules for classical groups. In
Section~\ref{sec:rank3} we collect the previous results to prove the vanishing
result in Theorem~\ref{thmA:zeroes}. Finally, in Section~\ref{thm:cartan} we
show the application to the proof of Theorem~\ref{thmA:PIMs}.
\medskip

\noindent
{\bf Acknowledgement:} We thank Andrew Mathas for useful explanations
concerning the notation and results of \cite{JM00}. We also thank the anonymous
referee for his helpful comments.

 \section{Maximal tori in classical groups}   \label{sec:maxtor}

In this section we establish the existence of regular semisimple $\ell$-singular
elements in suitable tori of classical groups of Lie type. This will be a
main ingredient in the proof of the vanishing Theorem~\ref{thmA:zeroes}.
Moreover, we formulate a criterion for collections of maximal tori of simple
algebraic groups not to be contained in a common proper reductive subgroup.

\subsection{Weyl groups and maximal tori in symplectic and orthogonal
groups}\label{subsec:WeylBCD}

For a prime power $q$ and $K = \overline{\FF_q}$ we consider the natural
representation of
$\bG = \Sp_{2n}(K)$ on $V = K^{2n}$ (case $C_n$) and of $\bG = \SO_{2n+1}(K)$
on $V = K^{2n+1}$, and a Steinberg morphism $F$ on $\bG$ with fixed
point group $\bG^F = \Sp_{2n}(q)$, respectively $\bG^F = \SO_{2n+1}(q)$. With
respect to an appropriate choice of basis of $V$ we can assume (see
\cite[Exmp.~6.7]{MT}) that the subgroup $\bT$ of diagonal matrices
forms an $F$-stable maximal torus of $\bG$, its elements have the form
$\diag(t_1, \ldots, t_n,t_n^{-1}, \ldots, t_1^{-1})$, respectively
$\diag(t_1, \ldots, t_n,1,t_n^{-1}, \ldots, t_1^{-1})$,
and $F$ acts by raising the diagonal entries to their $q$-th
power. The Weyl group $W$ acts on $\bT$ by its natural permutation action on
the diagonal entries. More precisely, $W$ is isomorphic to the wreath product
$C_2\wr\fS_n$, where the cyclic groups of order $2$ interchange the entries
$\{t_i, t_i^{-1}\}$, $i = 1, \ldots, n$, while the symmetric group
permutes the $n$ pairs $\{t_i, t_i^{-1}\}$. See~\cite[Chap.~15]{DM} for
more details.

The $\bG^F$-conjugacy classes of $F$-stable maximal tori of $\bG$ are
parameterised by the conjugacy classes of $W$. An element $t\in \bT$ is
conjugate to an element $t' \in \bG^F$ if ${}^wF(t)=t$ for some $w \in W$. Then
$t'$ lies in an $F$-stable maximal torus parameterised by the class of $w$.
Such an element $t'$ is regular semisimple if $w$ is unique.

The conjugacy classes of $W$ are parameterised by pairs of partitions
$(\lambda,\mu) \vdash n$ of $n$ as follows: If $w \in W$ permutes entries of
a regular diagonal element as above, then a cycle containing entries $t_i$ and
$t_i^{-1}$ for some $i$ is of even length $2\mu_j$ and contributes a part
$\mu_j$ to $\mu$. Otherwise a cycle of length $\lambda_j$ permutes some
diagonal entries and there is another cycle of the same length permuting the
inverse diagonal entries in the same way; this contributes a part
$\lambda_j$ to $\lambda$.

We can embed an even dimensional orthogonal group $\bG = \SO_{2n}(K)$
(type $D_n$) into $\tilde\bG = \SO_{2n+1}(K)$ such that it is stable under
the Steinberg morphism of $\tilde\bG$ and its fixed points are $\bG^F =
\SO^+_{2n}(q)$. There is a second Steinberg morphism $F'$ on $\bG$
leading to the twisted groups $\bG^{F'} = \SO^-_{2n}(q)$.
The Weyl group $W$ of $\bG$ is a subgroup of index $2$ in the Weyl group
$\tilde{W}$ of $\tilde\bG$, consisting of the conjugacy classes of $\tilde W$
whose parameter $(\lambda,\mu)$ has an even number of parts in $\mu$. If
$\mu$ is empty in $(\lambda,\mu)$  and $\lambda$ has only even parts then
the corresponding class in $\tilde W$ splits into two classes in $W$, called
\emph{degenerate} classes. The twisted Steinberg morphism  $F'$ acts on $W$
like conjugation by the short root generator of $\tilde W$. This way the
$F'$-conjugacy classes of $W$ can be parameterised by the conjugacy classes
of $\tilde W$ whose parameter $(\lambda, \mu)$ has an odd number of parts in
$\mu$. Since $\bG$ and $\tilde\bG$ have the same rank, maximal tori of $\bG$
are also maximal tori of $\tilde\bG$.

Let $(\la,\mu) = ((\la_1, \ldots, \la_r),(\mu_1,\ldots,\mu_s)) \vdash n$.
Then the order of the corresponding maximal torus $\bT^F$ in $\bG^F$ is
\[ |\bT^F| = \prod_{i=1}^{r}(q^{\la_i}-1) \prod_{i=1}^{s}(q^{\mu_i}+1) \]
and $\bT^F$ contains cyclic subgroups of orders $q^{\la_i}-1$ and
$q^{\mu_i}+1$ for all~$i$.

\subsection{Regular semisimple elements}\label{subsec:regss}
We need to know that certain maximal tori contain regular elements. As before,
let $q$ be a prime power. For a prime $\ell$ not dividing $q$ we denote by
$d_\ell(q)$ the order of $q$ modulo $\ell$, that is the smallest $d\ge1$ with
$\ell \mid (q^d-1)$. If $d$ is even then we have $\ell \mid (q^{d/2}+1)$.
With these definitions, $\ell \mid (q^k-1)$ if and only if $d \mid k$, and for
even $d$ we have $\ell \mid (q^k+1)$ if and only if $k$ is an odd multiple
of $d/2$.

\begin{lem}   \label{lem:regelts}
 Let $\bG$ be a simple simply-connected classical group of type $B_n$, $C_n$
 or $D_n$ defined over the finite field $\FF_q$ with corresponding Steinberg
 morphism $F$.

 Let $(\la,\mu) = ((\la_1, \ldots, \la_r),(\mu_1, \ldots, \mu_s))$ be a
 pair of partitions of $n$, and $\bT$ a corresponding
 $F$-stable maximal torus of $\bG$. Then $\bT^F$ contains regular elements
 if one of the following conditions is fulfilled.
 \begin{itemize}
  \item[(1)] $q>3$, $\la_1 < \la_2 < \ldots < \la_r$ and $\mu_1 < \mu_2 <
    \ldots < \mu_s$;
  \item[(2)] $q \in \{2,3\}$,
    $\la_1 < \la_2 < \ldots < \la_r$, $\mu_1 < \mu_2 <
     \ldots < \mu_s$, all $\la_i \neq 2$, and if $\bG$ is of type $B_n$ or
     $C_n$ then also all $\la_i \neq 1$; or
  \item[(3)] $\bG$ is of type $D_n$, $2 < \la_1 < \la_2 < \ldots < \la_r$ and
    $1 = \mu_1 = \mu_2 < \mu_3 < \ldots < \mu_s$.
 \end{itemize}
\end{lem}

\begin{proof}
(a) We use the natural representations of $\tilde\bG = \Sp_{2n}(K)$ and
$\tilde\bG = \SO_{2n+1}(K)$ with the natural Steinberg morphism that raises
matrix entries to their $q$-th power, and we consider $\SO_{2n}^\pm(q)$ as
subgroups of $\SO_{2n+1}(K)$.

We have described the maximal torus $\bT$ of diagonal elements in
$\tilde\bG$ in Section~\ref{subsec:WeylBCD}.
The connected centraliser $C_{\tilde\bG}^\circ(t)$ of $t\in\bT$ equals $\bT$
if and only if $\alpha(t) = 1$ for all roots $\alpha$ of $\bG$ with respect
to $\bT$. These roots are explicitly given in~\cite[Chap. 15]{DM} (for
types $C_n$ and $D_n$, but type $B_n$ is very similar): It is enough to
describe a set of positive roots. Their values on $t = \textrm{diag}(t_1,
t_2, \ldots, t_n, t_n^{-1}, \ldots, t_1^{-1})$ respectively
$t = \textrm{diag}(t_1, t_2, \ldots, t_n, 1, t_n^{-1}, \ldots, t_1^{-1})$ are

(Types $B_n$, $C_n$, $D_n$) $\alpha(t) = t_i t_j^{-1}$ and $\alpha(t) = t_i
t_j$ for $1\leq i < j \leq n$,

(Type $B_n$ only) $\alpha(t) = t_i$ for $1\leq i \leq n$,

(Type $C_n$ only) $\alpha(t) = t_i^2$ for $1\leq i \leq n$.

So, $C_{\tilde\bG}^\circ(t) = \bT$ if and only if for all $1 \leq i < j \leq n$
we have $t_i \neq t_j$ and $t_i \neq t_j^{-1}$, and in type $B_n$ we also
have $t_i \neq 1$ for $1 \leq i \leq n$, and in type $C_n$ we also have
$t_i \neq \pm 1$ for $1 \leq i \leq n$.

(b) In type $C_n$ we have $\bG = \tilde \bG$, while in types $B_n$ and $D_n$
the groups $\SO_{2n+1}(K)$ and $\SO_{2n}(K)$ are not simply-connected. But
there are isogenies from the simply-connected
coverings $\bG = \Spin_m(K) \to \SO_m(K) = \tilde\bG$. If $q$ is even, that
is $\textrm{char}(K) = 2$, then this map has a trivial kernel and induces
isomorphisms $\Spin^\pm_m(q) \to \SO_m^\pm(q)$. For odd $q$ the kernel is a
central subgroup of order $2$, and the image of the induced
map $\Spin^\pm_m(q) \to \SO_m^\pm(q)$ of the finite groups has index $2$.

To show the existence of regular elements in tori of type $(\la, \mu)$ we
construct a $t \in \bT$ such that $^wF(t) = t$ for some $w\in W$ in the class
$(\la,\mu)$ and $C_{\tilde\bG}^\circ(t) = \bT$. In cases $B_n$ and $D_n$ and odd
$q$ this $t$ will be a square of another element in $\bT^{wF}$ such that in
all cases $t$ has a preimage in a corresponding twisted maximal torus of the
simply-connected group $\bG$.  This preimage is regular because centralisers
of semisimple elements in a simply-connected group are connected.

(c) Now let $(\la,\mu) = ((\la_1, \ldots, \la_r),(\mu_1, \ldots, \mu_s))
\vdash n$ be a pair of partitions. We construct $t \in \bT$ by describing
the entries $t_1, \ldots, t_n \in K$.

First assume that $\bG$ is of type $C_n$ or that $q$ is even.

For $1 \leq i \leq r$ let $a_i \in \FF_{q^{\la_i}}$ be of order $q^{\la_i}-1$
and use $a_i, a_i^q, \ldots, a_i^{q^{\la_i-1}}$ as entries of $t$ (note that
$a_i^{q^{\la_i}} = a_i$).

For $1 \leq j \leq s$ let $b_j \in \FF_{q^{2\mu_j}}$ be of order
$q^{\mu_j}+1$ and use $b_j, b_j^q, \ldots, b_j^{q^{\mu_j-1}}$ as entries of
$t$ (note that then $b_j^{q^{\mu_j}} = b_j^{-1}$).

If $\bG$ is of type $B_n$ or $D_n$ and $q$ is odd, we use as $t$ the square
of the element just described (so $a_i$ is of order $\frac{1}{2}(q^{\la_i}-1)$
and $b_j$ is of order $\frac{1}{2}(q^{\mu_j}+1)$).

It is clear from this construction that the Steinberg morphism $F$ of
$\Sp_{2n}(K)$, respectively $\SO_{2n+1}(K)$ permutes the entries of $t$ such
that this can be reversed by an element $w \in W$ of cycle type $(\la,\mu)$,
so $t = {}^wF(t)$.

(d) We next discuss when $t$ fails to be regular. Let $\la_i,\mu_j>0$ with
$q^{\la_i}-1 = q^{\mu_j}+1$. Then $\mu_j < \la_i$ and the assumption is
equivalent to $q^{\mu_j}(q^{\la_i-\mu_j}-1) = 2$. It follows that $q=2$,
$\mu_j=1$ and  $\la_i=2$.

(e) Let $k \in \ZZ$, $0 < k < \la_i$ and assume $(q^{\la_i}-1)\mid(q^k \mp 1)$.
Then $q^{\la_i}-1 \leq q^k \mp 1$, which holds
if and only if $q^k(q^{\la_i-k} -1) \leq 0$ or $2$, respectively. This is only
possible for $q=2$, $k=1$ and $\la_i=2$. Similarly, for odd $q$ the condition
$\frac{1}{2}(q^{\la_i}-1)\mid (q^k \mp 1)$ implies $q=3$, $k=1$ and $\la_i=2$.

Now let $k \in \ZZ$, $0 < k < \mu_j$.
A similar argument as above shows that for all $q$ we have
$(q^{\mu_j}+1) \nmid (q^k \mp 1)$ and for odd $q$ we also have
$\frac{1}{2}(q^{\mu_j}+1) \nmid (q^k \mp 1)$.

(f) Now we prove assertions~(1) and~(2) of the statement by showing that in
those cases the element $t\in\bT$ constructed in~(c) does not lie in the kernel
of any root.

First note that the $a_i$ ($1\leq i \leq r$) and $b_j$ ($1 \leq j \leq s$)
have pairwise different orders. This follows from~(d) since for $q=2$ we
assume that all $\la_i \neq 2$.

According to~(a)  we need to show that for $1\leq i \leq r$ any of the entries
$(a_i, a_i^q, \ldots, a_i^{q^{\la_i-1}})$ are not equal or inverse to each
other, and similarly for the $(b_j, b_j^q, \ldots, b_j^{q^{\mu_j-1}})$.
Since the map $x \mapsto x^q$ permutes these entries of $t$ and their
inverses it suffices to show that for $0 <k <\la_i$ we have $a_i^{q^k} \neq
a_i^{\pm 1}$ (equivalently $a_i^{q^k\mp 1} \neq 1$), and similarly for
the $b_j$. Now $a_i^{q^k\mp 1} = 1$ implies $(q^{\la_i}-1) \mid (q^k \mp
1)$ or $\frac{1}{2}(q^{\la_i}-1) \mid (q^k \mp 1)$, respectively. Using~(e)
this implies $\la_i = 2$ and $q\in \{2,3\}$, contradicting the assumptions.

Similarly, step~(e) shows that we always have $b_j^{q^k} \neq b_j^{\pm 1}$
for $0 < k < \mu_j$.

For type $C_n$ we also have to show that $a_i^2 \neq 1$ ($1\leq i \leq r$)
and $b_j^2 \neq 1$ ($1\leq j\leq s$). This holds because $a_i^2 =1$,
i.e., $(q^{\la_i} -1) \mid 2$, implies $\la_i =1$ and $q \in
\{2,3\}$, contradicting the assumption. Note that always $q^{\mu_i}+1 > 2$.

For type $B_n$ we also have to show that $a_i \neq 1$ ($1\leq i \leq r$)
and $b_j \neq 1$ ($1\leq j\leq s$). This follows as $q^{\la_i}-1 = 1$ implies
$q=2$ and $\la_i=1$, while $\frac{1}{2}(q^{\la_i}-1) = 1$ implies $q=3$ and
$\la_i=1$ for odd $q$, and as always $q^{\mu_j}+1 >1$.

(g) The argument for statement~(3) is very similar. A suitable element $t$
can be defined as before, except that for $\mu_1=1$ we choose $b_1  = 1$.
Then the entry for $\mu_2=1$ has order $q+1$ for even $q$ or
$\frac{1}{2}(q+1)$ for odd $q$ which is always larger than $1$.
\end{proof}

\begin{rem}   \label{rem:lsingreg}
Let $\bG$ be a connected reductive group with Steinberg morphism $F$ and
$\bT$ an $F$-stable maximal torus of $\bG$. Let $\ell$ be a prime dividing
$|\bT^F|$. If $\bT^F$ contains a regular element then $\bT^F$ also contains
an $\ell$-singular regular element.
\end{rem}

\begin{proof}
Let $t\in \bT^F$ be an $\ell$-regular (that is of order prime to $\ell$)
regular element  and $u\in \bT^F$ be any element of order $\ell$. Then $tu$
has order $\ell |t|$ and some power of $tu$ equals $t$, hence $tu$ is also
regular.
\end{proof}

\subsection{Maximal connected reductive subgroups}\label{subsec:maxred}
The following lemma will allow us to decide if a collection of classes of
maximal tori of a classical group has representatives inside the centraliser
of a semisimple element:

\begin{lem}   \label{lem:maxred}
 Let $\bG$ be a simple algebraic group of classical type $B_n,C_n$ (with
 $n\ge2$) or $D_n$ (with $n\ge4$) with Frobenius endomorphism $F$ such that
 $\bG^F$ is a classical group. Let $\Lambda$ be a set of pairs of partitions
 $(\la,\mu)\vdash n$. Assume the following:
 \begin{enumerate}
  \item[\rm(1)] there is no $1\le k\le n-1$ such that all $(\la,\mu)\in\Lambda$
   are of the form $(\la_1,\mu_1)\sqcup(\la_2,\mu_2)$ (that is $\la$
   contains the parts from $\la_1$ and $\la_2$ and similarly for $\mu$) with
   $(\la_1,\mu_1)\vdash k$;
  \item[\rm(2)] the greatest common divisor of all parts of all
   $(\la,\mu)\in\Lambda$ is~1; and
  \item[\rm(3)] if $\bG$ is of type $B_n$, then there exist pairs
   $(\la,\mu)\in\Lambda$ for which $\mu$ has an odd number of parts, and for
   which $\mu$ has an even number of parts.
 \end{enumerate}
 If $s\in\bG^F$ is semisimple such that $C_\bG(s)$ contains maximal tori of
 $\bG$ corresponding to all $(\la,\mu)\in\Lambda$ then $s$ is central.
\end{lem}

\begin{proof}
The types of maximal tori are insensitive to isogeny types, so we may assume
that $\bG$ is a classical matrix group $\Sp(V)$ or $\SO(V)$. Let
$s\in\bG^F$ be a non-central semisimple element.
Then $\bH:=C_\bG^\circ(s)$ fixes each eigenspace of $s$ on $V$. Since $s$ is
non-central it has at least two different eigenvalues. Since $s$ is $F$-stable,
the eigenspaces of $s$ are permuted under the action of $F$. First
assume that the set of eigenspaces forms one orbit, of length~$f>1$ say, under
the action of $F$. Then $\bH$ is a central product of isomorphic classical
groups, $\bH^F$ is an extension field subgroup and hence its maximal tori
have parameters all parts of which are divisible by $f$, contradicting~(2).
\par
So there are at least two $F$-orbits of eigenspaces. Let $V_1< V$ be an
$F$-stable sum of eigenspaces of minimal possible dimension. Then either
$V_1 \cap V_1^\perp \neq 0$ or $V_1$ is non-degenerate and $V = V_1 \perp
V_1^\perp$. The stabilisers
of totally singular spaces are contained in maximal parabolic subgroups
\cite[Prop.~12.13]{MT}, but these only have tori of types as excluded by~(1).
Else, $\bH$ is contained in the stabiliser of $V_1\perp V_1^\perp$, a central
product of
classical groups. The latter are of types $C_a+C_b$, respectively $B_a + D_b$
or $D_a + D_b$ with $a+b = n$, $a,b > 0$ for $\bG$ of type $C_n$, respectively
$B_n$, $D_n$. Then the parameters of $F$-stable maximal tori of $\bH$ are of
the form excluded by~(1) of the lemma, for $k=a$, unless $\bG$ is of type
$B_n$ and $V_1$ is 1-dimensional, when $\bH$ could be of type $D_n$.  In
that case, the number of parts in the second entry of the parameters of
$F$-stable maximal tori of $\bH$ all have a fixed parity (namely even for
$\bH^F$ untwisted, and odd for $\bH^F$ twisted). This is excluded by~(3).
\end{proof}

\section{A Murnaghan--Nakayama rule}   \label{sec:MNrule}

In this section we derive a Murnaghan--Nakayama type formula for the values
of unipotent characters of classical groups of Lie type on regular semisimple
elements. It relies on a result of Asai on the decomposition of Lusztig
restriction.

\subsection{Lusztig restriction}  \label{subsec:*RLG}
We consider the following setup: $\bG$ is a connected reductive algebraic
group over an algebraic closure of the finite field $\FF_p$, and
$F:\bG\rightarrow\bG$ is a Steinberg endomorphism of $\bG$, with finite
group of fixed points $\bG^F$. For an $F$-stable Levi subgroup $\bL\le\bG$
of $\bG$ Lusztig defined, via certain $\ell$-adic cohomology constructions,
induction and restriction functors
$$\RLG:\ZZ\Irr(\bL^F)\rightarrow\ZZ\Irr(\bG^F)\quad\text{ and }\quad
  \sRLG:\ZZ\Irr(\bG^F)\rightarrow\ZZ\Irr(\bL^F)$$
between the respective groups of virtual characters. These are adjoint to one
another with respect to the usual scalar product of characters. We need the
following connection between Lusztig restriction and ordinary restriction on
regular semisimple conjugacy classes:

\begin{prop}   \label{prop:restr}
 Let $\bL\le\bG$ be an $F$-stable Levi subgroup. Let $\chi\in\Irr(\bG^F)$ and
 $s\in \bL^F$ be regular semisimple. Then
 $\sRLG(\chi)(s)=\chi(s)$.
\end{prop}

\begin{proof}
According to the character formula in \cite[Prop.~12.2(ii)]{DM} we have
$$\sRLG(\chi)(s)=\frac{|C_\bL^\circ(s)^F|}{|C_\bG^\circ(s)^F|}
   \sum_{u\in C_\bG^\circ(s)^F_u}Q_{C_\bL^\circ(s)}^{C_\bG^\circ(s)}(u,1)
   \,\chi(s),$$
where the sum runs over unipotent elements in $C_\bG^\circ(s)^F$, and
$Q_{C_\bL^\circ(s)}^{C_\bG^\circ(s)}$ is a certain 2-parameter Green function.
Since $s$ is assumed to be regular semisimple,
$\bT:=C_\bG^\circ(s)=C_\bL^\circ(s)$ is a torus. Hence the only term in the
sum is the one for $u=1$, and by \cite[p.~98]{DM} the Green function takes
the value
$$Q_{C_\bL^\circ(s)}^{C_\bG^\circ(s)}(1,1)=Q_{\bT}^{\bT}(1,1) =1$$
(again using that $C_\bG^\circ(s)$ contains no non-trivial unipotent elements).
The claim follows.
\end{proof}

\subsection{Symbols}  \label{subsec:symb}
We now specialise to simple groups of classical type. By Lusztig's results
their unipotent characters are parameterised by so-called \emph{symbols}.
Here, a symbol is an unordered pair of finite sets of non-negative integers,
usually denoted as
$$\cS=\binom{X}{Y}=\begin{pmatrix} \la_1& \la_2& \ldots& \la_r\\
  \mu_1& \mu_2& \ldots& \mu_s\end{pmatrix}$$
with $\la_1<\la_2<\ldots<\la_r$ and $\mu_1<\mu_2<\ldots<\mu_s$ where $r,s\ge0$.
Two symbols are said to be \emph{equivalent} if one can be obtained from the
other by a sequence of shift operations
$$\cS\mapsto \begin{pmatrix} 0& \la_1+1& \la_2+1& \ldots& \la_r+1\\
  0& \mu_1+1& \mu_2+1& \ldots& \mu_s+1\end{pmatrix}$$
or by interchanging the rows $X$ and $Y$ of $\cS$. The \emph{rank} of $\cS$ is
by definition
$$\rk(\cS)=\sum_{i=1}^r \la_i+\sum_{i=1}^s\mu_i
  -\Big\lfloor\left(\frac{r+s-1}{2}\right)^2\Big\rfloor.$$
The \emph{defect} of $\cS$ is $\df(\cS):=|r-s|$. Note that the rank and the
defect are well-defined on equivalence classes.
We will use the following notation to describe explicit symbols:
$\ul{n} := \{0,1,\ldots,n\}$ and for a set $X$, a $k\in X$ and $l\notin X$ we
write $X \sm k \cup l := X \sm\{k\} \cup \{l\}$.

\par
Let $d$ be a positive integer. We say that $\cS$ has a \emph{$d$-hook} $h$ at
$x\in X$ if $0\le x-d\notin X$, and similarly at $x\in Y$. Removing that
$d$-hook leads to the symbol $\cS\sm{h}:=\binom{X'}{Y}$ with
$X'=X\sm x\cup(x-d)$ (respectively $\cS\sm{h}:=\binom{X}{Y'}$
with $Y'=Y\sm x\cup(x-d)$). Attached to the hook $h$ is the sign
$\eps_h:=(-1)^m$ where $m:=|\{y\in X\mid x-d<y<x\}|$ (respectively
$m:=|\{y\in Y\mid x-d<y<x\}|$). We say that $\cS$ has a \emph{$d$-cohook} $c$
at $x\in X$ if $0\le x-d\notin Y$, and similarly at $x\in Y$.  Removing that
$d$-cohook leads to the symbol $\cS\sm{c}:=\binom{X'}{Y'}$ with
$X'=X\sm\{x\},Y'=Y\cup\{x-d\}$ (respectively $Y'=Y\sm\{x\},
X'=X\cup\{x-d\}$). Attached to the cohook $c$ is the sign $\eps_c:=(-1)^m$,
where $m:=|\{y\in X\mid y<x\}|+|\{y\in Y\mid y<x-d\}|$ (respectively
$m:=|\{y\in Y\mid y<x\}|+|\{y\in X\mid y<x-d\}|$), see \cite[\S3]{FS86} for
example. Note that removing a $d$-hook or a $d$-cohook from a symbol of rank
$n$ yields a symbol of rank $n-d$.

For a symbol $\cS=(X,Y)$ with all entries of $X,Y$ not larger than $m$ we
define the \emph{dual symbol} $\cS^\vee=(X',Y')$ where
$$X'= \ul{m} \setminus \{m-x\mid x\in X\}\mbox{ and }
  Y'= \ul{m} \setminus \{m-y\mid y\in Y\}.$$

The hooks and cohooks of $S$ and $S^\vee$ are each in natural bijection such
that removing hooks or cohooks commutes with the duality operation.
This observation will allow us to simplify some of the later statements and
proofs.

\subsection{Unipotent characters of symplectic and orthogonal groups}   \label{subsec:uniBCD}
Now assume that $\bG$ is simple of type $B_n$, $C_n$ or $D_n$, and $F$ does
not induce the triality automorphism when $\bG$ is of type $D_4$ nor the
exceptional graph automorphism in type $B_2$ or $C_2$. According to Lusztig
\cite{Lu84} the unipotent characters of $G:=\bG^F$ are then parameterised by
(equivalence classes of) symbols of rank~$n$. More precisely, if $\bG$ is of
type $B_n$ or $C_n$, the unipotent characters of $G$ are parameterised by
symbols of rank~$n$ and odd defect. Now assume that $\bG$ is of type $D_n$.
If the Steinberg endomorphism $F$ is untwisted, then the unipotent characters
of $G$ are parameterised by symbols $\cS$ of rank~$n$ and defect
$\df(\cS)\equiv0\pmod4$, where symbols with two equal rows stand for two
characters each. These are the so-called \emph{degenerate} symbols. If $\bG^F$
is of twisted type, then the unipotent characters of $G$ are parameterised
by symbols $\cS$ of rank~$n$ and defect $\df(\cS)\equiv2\pmod4$.

If $\cS$ is a symbol parameterising a unipotent character $\rho$ of a classical
group then the Alvis--Curtis dual $\rho^\vee$ of $\rho$ is parameterised by the
dual symbol $\cS^\vee$ which has the same rank, the same defect and the same
hook and cohook lengths (see \cite[4.5.5, 4.6.8]{Lu84}).
(In fact, the degree $\rho^\vee(1)$ differs by a power of the underlying
characteristic of $\bG$ from the degree $\rho(1)$.)

\subsection{Values of unipotent characters at regular semisimple
elements}\label{subsec:MN}
We first recall the following result of Asai on the decomposition of Lusztig
induction as given in~\cite[(3.1), (3.2)]{FS86}. We restate it in terms of the
restriction functors $\sRLG$ instead of the induction functors $\RLG$ as
in the cited reference. This is a summary of results from~\cite[2.8]{As1},
\cite[1.5]{As2} and~\cite[2.2.3]{As3}. Note that there is a sign missing
in the formula~\cite[(3.2)]{FS86} in the case $D_n$ as can be seen by
considering the trivial character.

\begin{thm}[Asai (1984)]   \label{thm:Asai}
 Let $\bG$ be of type $B_n,C_n$ or $D_n$. Let $d\ge1$.
 Let $\cS$ be a symbol of a unipotent character of $G=\bG^F$. If $\cS$ is not
 degenerate, we write $\rho_\cS$ for the corresponding character. If $\cS$
 is degenerate we write $\rho_\cS$ for the sum of the two corresponding
 characters.
 \begin{enumerate}[\rm(a)]
  \item Let $\bL$ be the centraliser in $\bG$ of an $F$-stable torus $\bT_d$
   of $\bG$ with $|\bT_d^F|=q^d-1$. Then
   $$\sRLG(\rho_\cS)=\sum_{h\text{ $d$-hook}}\eps_h\,\rho_{\cS\setminus h}$$
   where $h$ runs over the $d$-hooks of $\cS$.
  \item Let $\bL$ be the centraliser in $\bG$ of an $F$-stable torus $\bT_d$
   of $\bG$ with $|\bT_d^F|=q^d+1$. Then
   $$\sRLG(\rho_\cS)=(-1)^\delta\sum_{c\text{ $d$-cohook}}
     \eps_c\,\rho_{\cS\setminus c}$$
   where $c$ runs over the $d$-cohooks of $\cS$, and $\delta=0$ for types
   $B_n,C_n$, $\delta=1$ for type $D_n$.
 \end{enumerate}
\end{thm}

With this we can show a Murnaghan--Nakayama formula for values of unipotent
characters at regular semisimple elements.

\begin{thm}   \label{thm:MN}
 Let $\bG$ be of type $B_n,C_n$ or $D_n$.
 Let $\cS$ be a symbol of a unipotent character of $G=\bG^F$. If $\cS$ is not
 degenerate, we write $\rho_\cS$ for the corresponding character. If $\cS$
 is degenerate we write  $\rho_\cS$ for the sum of the two corresponding
 characters.
 Let $s\in G$ be regular semisimple lying in an $F$-stable maximal torus $\bT$
 of $\bG$ parameterised by $(\la,\mu)\vdash n$.
 \begin{enumerate}[\rm(a)]
  \item If $\la$ has a part of length $d$ then
   $$\rho_\cS(s)=\sum_{h\text{ $d$-hook}}\eps_h\,\rho_{\cS\setminus h}(s)$$
   where $h$ runs over the $d$-hooks of $\cS$.
  \item If $\mu$ has a part of length $d$ then
   $$\rho_\cS(s)=(-1)^\delta\sum_{c\text{ $d$-cohook}}
     \eps_c\,\rho_{\cS\setminus c}(s)$$
   where $c$ runs over the $d$-cohooks of $\cS$, and $\delta=0$ for types
   $B_n,C_n$, $\delta=1$ for type $D_n$.
 \end{enumerate}
 On the right hand side of these formulae the symbols and characters belong
 to a Levi subgroup of semisimple rank $n-d$ and $s$ lies in a maximal torus
 of that Levi subgroup corresponding to the pair of partitions
 $(\lambda \setminus d,\mu) \vdash (n-d)$ in case~(a), respectively
 $(\lambda, \mu\setminus d) \vdash (n-d)$ in case~(b).
\end{thm}

\begin{proof}
First assume that $\la$ has a part of length $d$. Then the torus $\bT$ has
an $F$-stable subtorus $\bT_d$ with $|\bT_d^F|=q^d-1$, see
Section~\ref{subsec:WeylBCD}. We consider Lusztig restriction of
$\rho=\rho_\cS$ to the $F$-stable Levi subgroup
$\bL:=C_\bG(\bT_d)\cong\bH\bT_d$, with $\bH$ a group of the same type as
$\bG$ and of rank~$n-d$. This is a $d$-split Levi subgroup containing $\bT$,
so $s$, and according to Asai's Theorem~\ref{thm:Asai}(a), the
constituents of $\sRLG(\rho)$
are, up to the sign $\eps_h$ as given in the statement, precisely those
unipotent characters of $\bL^F$ whose parameterising symbol is obtained from
$\cS$ by removing a $d$-hook. Application of Proposition~\ref{prop:restr}
then gives the claim in~(a).
\par
If $\mu$ has a part of length $d$, then $\bT$ has an $F$-stable subtorus
$\bT_d$ with $|\bT_d^F|=q^d+1$, and we can argue precisely as before using
Theorem~\ref{thm:Asai}(b) for the decomposition of Lusztig restriction to
$\bL:=C_\bG(\bT_d)\ge\bT\ni s$.
\end{proof}

Note that in the theorem $\rho_{\cS\sm h}$, respectively $\rho_{\cS\sm c}$,
is a unipotent character of the classical group $\bL^F$ of rank $n-d$ (of the
same type as $G$, except that removing cohooks changes the value of the defect
modulo~4, so in type $D_n$ interchanges the twisted and untwisted types), and
that $s$ is a fortiori regular in $\bL$. Thus, the above result gives a
recursive algorithm to compute the character values on all regular semisimple
elements.

In case of degenerate symbols we only get the values of the sum of two
unipotent characters. These unipotent characters form one element Lusztig
families and so their values on regular semisimple classes are the values of
the corresponding degenerate irreducible characters of the Weyl group,
see~\cite[(4.6.10)]{Lu84}. The two characters for
a degenerate symbol have the same values on non-degenerate classes, but
different values on the degenerate classes. See~\cite{PfWeyl} for a
definition of unique parameters for the degenerate characters and classes
and a formula to compute the missing values on degenerate classes.

The following vanishing result is immediate from the above:

\begin{cor}   \label{cor:vanish}
 Let $\bG$, $\cS$ and $\rho_\cS$ be as in Theorem~\ref{thm:MN}, and
 let $s\in G$ be a regular semisimple element lying in an $F$-stable
 maximal torus parameterised by $(\la,\mu)$.
 Then $\rho_\cS(s)=0$ in any of the two following cases:
 \begin{enumerate}[\rm(1)]
  \item $\la$ has a part of length $d$, but $\cS$ does not have a $d$-hook; or
  \item $\mu$ has a part of length $d$, but $\cS$ does not have a $d$-cohook.
 \end{enumerate}
 Thus, if $\rho_\cS(s) \neq 0$ and $(\la,\mu) =
 ((\la_1, \ldots, \la_r),(\mu_1,\ldots,\mu_s))$, then it is possible to
 remove $\la_i$-hooks for $1\leq i \leq r$ and $\mu_i$-cohooks for
 $1\leq i \leq s$ from $\cS$ in any chosen order.
\end{cor}

We remark that we have computed many character values on semisimple classes
for classical groups of rank $\leq 10$ by Deligne--Lusztig theory. These
examples helped to find the statements given in the next section. And they
also provided an independent check of the Murnaghan--Nakayama
formula~\ref{thm:MN}.

\section{Vanishing of characters}   \label{sec:vanish}
In this section we prove Theorem~\ref{thmA:zeroes} for quasi-simple groups
of classical type. So let $\bG$ be a simply-connected simple algebraic group of
type $B_n$ with $n \geq 3$, $C_n$ with $n \geq 2$ or $D_n$ with $n \geq 4$,
and let $F$ be a Steinberg morphism of $\bG$ (we exclude the $\tw3D_4$ case)
with $G =\bG(q) := \bG^F$. We consider the corresponding finite groups
$G=\Spin_{2n+1}(q)$ only for odd $q$, and $\Sp_{2n}(q)$ and $\Spin^\pm_{2n}(q)$
for any $q$, respectively.

We consider primes $\ell$ such that the Sylow $\ell$-subgroups of $G$ are not
cyclic. We will show that for such $\ell$ with only very few exceptions all
ordinary irreducible characters of $G$ vanish on some $\ell$-singular regular
semisimple element.

If we write $d := d_\ell(q)$ for the order of $q$ modulo $\ell$ as in
Section~\ref{subsec:regss} and then $n = ad+r$, $0\leq r < d$, for odd $d$
respectively $n=ae+r$, $0\leq r < e$, for even $d=2e$, then our assumption on
$\ell$ implies $a\geq 2$ because the cyclotomic polynomial $\Phi_d(q)$ must
divide the order of $\bG(q)$ at least twice.
If $G$ is of type $\tw2D_n$ we have $n > 2d$ for odd $d$ and $n > 2e$ for even
$d=2e$ (see the order formula in~\cite[Tab.~24.1]{MT}).

\subsection{Non-unipotent characters}   \label{ssec:nonunipotent}

We first deal with non-unipotent characters, and afterwards examine in
more detail the unipotent characters.

\begin{thm}   \label{thm:non-unip}
 Let $G$ be one of the groups $\Spin_{2n+1}(q)$ for odd $q$ and $n \geq 3$,
 $\Sp_{2n}(q)$ for any $q$ and $n \geq 2$, or $\Spin^\pm_{2n}(q)$ for any $q$
 and $n \geq 4$. Let $2\ne\ell{\not|}q$ be a prime such that the Sylow
 $\ell$-subgroups of $G$ are non-cyclic. Then any non-unipotent irreducible
 complex character of $G$ vanishes on some $\ell$-singular regular semisimple
 element, except for the two cases $G=C_2(2)=\Sp_4(2)$ and $G=C_4(2)=\Sp_8(2)$.

 The group $\Sp_4(2)\cong\fS_6$ contains no $3$-singular regular semisimple
 elements. It has three non-trivial irreducible characters which do not vanish
 on any $3$-singular class, the unipotent character $\sys{0 1 2}{-}$ and the
 two characters of degree $10$ (belonging to Lusztig series of type
 $A_1(q)\times (q+1)$).

 The group $\Sp_8(2)$ has only one $5$-singular regular semisimple class and
 there are $15$ irreducible characters which do not vanish on this class.
 There is only one non-trivial character which does not vanish on any
 $5$-singular class, the unipotent character $\sys{0 1}{4}$.
\end{thm}

\begin{proof}
We use the following facts from Deligne--Lusztig theory. The irreducible
characters of $G$ are partitioned into Lusztig series which are parameterised
by semisimple conjugacy classes of the dual group $G^*$.
If $s \in G^*$ is semisimple then the values of characters in the Lusztig
series of $s$ on semisimple elements $t \in G$ are linear combinations of
values of Deligne--Lusztig characters on $t$ which belong to types of tori
which occur in the connected centraliser of $s$ in $\bG^*$. In particular
Deligne--Lusztig characters corresponding to a maximal torus $T$ of $G$
vanish on regular semisimple classes which do not intersect $T$, see e.g.
\cite[Prop.~6.4]{LMS13}.

Since $\bG$ is simply-connected the dual group $\bG^*$ has trivial center
so that non-unipotent characters are in Lusztig series of non-central
elements.

We show the theorem by listing in each case a set of types of maximal tori
which contain regular $\ell$-singular elements  by Lemma~\ref{lem:regelts} and
which cannot all occur in any proper centraliser. We have described
these subgroups together with the types of tori they contain in
Lemma~\ref{lem:maxred}. Note that for $\bG$ of type $B_n$ the dual group is
of type $C_n$ (we need to apply Lemma~\ref{lem:maxred} for type $C_n$ and
consider tori with regular elements in type $B_n$); the same remark holds
with $B_n$ and $C_n$ interchanged.

Let $d=d_\ell(q)$ and write $n = ad+r$ with $0 \leq r < d$ for odd $d$ and
$n = ae+r$ with $0 \leq r < e$ for even $d=2e$, $a\geq 2$, as explained above.

Recall that for odd $d$ a torus of type $(\la,\mu)$ contains $\ell$-singular
elements if and only if some part of $\la$ is a positive multiple of $d$.
For even $d=2e$ a torus of type $(\la,\mu)$ contains $\ell$-singular elements
if and only if $\mu$ contains an odd multiple of $e$ or $\la$ contains an
even multiple of~$e$.

We collect the tori we consider in Table~\ref{tab:non-unip}.

\begin{table}\label{tab:non-unip}
\caption{Tori needed in the proof of Theorem~\ref{thm:non-unip}}
\begin{tabular}{|l|c|c|c|p{0.63\textwidth}|}
\hline
type&$d$&$a$&$r$&tori\\
\hline
$C_n$,$B_n$&odd&any&$0$&$((n),-)$, $((n-d),(d))$, and $((n-d),(1,d-1))$ if
$d>1$\\
$C_n$,$B_n$&odd&any&$>0$&$((n-r),(r))$, $((n-r-d),(d+r))$,
$((n-r-d),(d+r-1, 1))$\\
$D_n$&odd&any&0&$((n),-)$, $((n-1,1),-)$ if $d=1$, or $((n-d),(d-1,1))$ if
$d>1$\\
$D_n$&odd&any&$>0$&$((n-d-r, d+r),-)$, $((n-r,r),-)$ if $r \neq 2$, or
$((n-2),(1,1))$ if $r=2$\\
$^2D_n$$(^*)$&odd&any&0&$((n-d),(d))$, $((d),(n-d))$, $((n-d-1,d),(1))$ if
$d>1$\\
$^2D_n$&odd&any&$>0$&$((n-r),(r))$, $((n-d-r),(d+r))$, $((n-d-r,1),(d+r-1))$\\
$C_n$,$B_n$&$=2e$&even&$0$&$((n),-)$, $((n-e),(e))$, $(-,(n-e-1,e,1))$ if
$e>1$\\
$C_n$,$B_n$&$=2e$&even&$>0$&$((n-r),(r))$, $(-,(n-e,e))$, $(-,(n-e-1,e,1))$
if $r\neq 1$\\
$C_n$,$B_n$&$=2e$&odd&$0$&$(-,(n))$, $(-,(n-e,e))$, $((n-e-1),(e,1))$ if
$e>1$\\
$C_n$,$B_n$&$=2e$&odd&$>0$&$(-,(n-r,r))$, $((n-r-e),(r+e))$,
$((n-r-e),(r+e-1,1))$\\
$D_n$&$=2e$&even&$0$&$((n),-)$, $((n-e-1),(e,1))$ if $(n,e) \notin
\{(4,1),(6,3)\}$, or $((n-e-2),(e,2))$ otherwise\\
$D_n$&$=2e$&even&$>0$&$(-,(n-e,e))$, $((n-1,1),-)$ if $r=1$,
$((n-r),(r-1,1))$ if $r>1$\\
$D_n$$(^*)$&$=2e$&odd&$0$&$(-,(n-e,e))$, $(-,(n-2e,2e))$,
$((1), (n-2e, 2e-1))$ if $e>1$, $((2e),(e-1,1)$ if $n=3e$ \\
$D_n$&$=2e$&odd&$>0$&$(-,(n-r,r))$, $(-,(n-e,e))$, $((1),(n-r,r-1))$ if
$r>1$\\
$^2D_n$&$=2e$&even&$0$&$((n-e),(e))$, $((n-2e),(2e))$, $((n-2e, 2e-1),(1))$
if $e>1$\\
$^2D_n$&$=2e$&even&$>0$&$((n-r),(r))$, $((n-e),(e))$, $((n-r,1),(r-1))$ if
$r>1$\\
$^2D_n$&$=2e$&odd&$0$&$(-,(n))$, $((n-e),(e))$, $((n-e,1),(e-1))$ if $e>1$\\
$^2D_n$&$=2e$&odd&$>0$&$((n-e-r),(e+r))$, $(-,(n-r,r-1,1))$
if $r>1$, $((1),(n-1))$ if $r=1$\\
\hline
\end{tabular}
\end{table}

A few additional remarks are in order. In the two cases marked with $(^*)$ in
the table we need to rule out the possibility that the given tori all lie in
a subgroup of type $D_d + D_{n-d}$ or $D_e + D_{n-e}$, respectively.
Such a subgroup stabilises a decomposition of the orthogonal space on which $G$
acts into two non-degenerate subspaces. In the
first case one of these subspaces must be of minus- and the other of plus-type.
So the corresponding finite subgroups are of type $\tw2D_d+D_{2d}$ or
$D_d+\tw2D_{2d}$, each of these only contains one of the types of tori,
$((d),(n-d))$ or $((n-d),(d))$. The argument for the second case is similar,
there both subspaces must be of minus-type or both of plus-type.

It remains to check for which groups one of the tori given in
Table~\ref{tab:non-unip} does not contain regular elements.
Note that for $q=3$ we have $d > 2$ because $\ell \neq 2$, and for $q=2$ we
have $d > 1$. Furthermore, when $r>0$, we have $d\geq 3$ or $e \geq 2$,
respectively. This implies that only the following cases need an extra
consideration:

$C_2(2)$ with $e=1$: this is the first exception mentioned in the statement.

$B_4(q)$ with $q=3$ and $e=2$: We compute that for $q=3$ the torus of type
$((1),(1,2))$ contains regular element. So, we can argue with tori of types
$((4),-)$ and $((1),(1,2))$.

$C_4(q)$ with $q=3$ and $e=2$: In this case the torus of type $((2),(2))$
in $C_4(3)$ contains regular elements.
Here we can argue with tori of types $((4),-)$
and $((2),(2))$ by checking that the dual group (adjoint of type $B_4(3)$)
does not contain a semisimple element whose centraliser contains tori of
both types.

For $C_4(2)$ and $e=2$ (and so $\ell = 5$) only tori of type
$((4),-)$ contain $\ell$-singular regular elements.
This leads to the second exception mentioned in the statement.

$\tw2D_4(2)$ with $e=1$: Here, the maximal tori containing $\ell$-elements all
have conjugates inside a subgroup of type $D_1+D_3$. We can apply the argument
from above: the tori of types $((1),(3))$ and $((3),(1))$ cannot both occur in
the same rational form of such a subgroup.
\end{proof}

\subsection{Unipotent characters}   \label{ssec:unipotent}
In the following theorem we use symbols to describe unipotent characters
and refer to Section~\ref{subsec:symb} for the notation.

\begin{thm}   \label{thm:unip}
 Let $G$ be one of the groups $\Spin_{2n+1}(q)$ for odd $q$ and $n \geq 3$,
 $\Sp_{2n}(q)$ for any $q$ and $n \geq 2$, or $\Spin^\pm_{2n}(q)$ for any $q$
 and $n \geq 4$. Let $2\ne\ell{\not|}q$ be a prime such that the Sylow
 $\ell$-subgroups of $G$ are non-cyclic.

 Let $\rho_\cS$ be a unipotent character of $G$ for some symbol $\cS$ and
 assume that $\rho_\cS$ does not vanish on any regular semisimple
 $\ell$-singular class. Then either $\rho_\cS$ is the trivial character, or
 the Steinberg character, or up to Alvis-Curtis duality we are in one of the
 following cases (as before we set $d = d_\ell(q)$, and interpret $\ul{-1}$ as
 the empty set):
 \begin{enumerate}[\rm(1)]
  \item $G=\Sp_{2n}(q)$ or $\Spin_{2n+1}(q)$, $d$ is odd, $n = 2d+r$ with
   $0\leq r < d$ and
   \[\cS = \sy{\ul{d-r-1}\sm 0 \cup d \cup 2d}{\ul{d-r-1}};\]
  \item $G=\Sp_{2n}(q)$ or $\Spin_{2n+1}(q)$, $d=2e$ is even, $n=2e+r$ with
   $0\leq r < e$ and
   \[\cS = \sy{\ul{e-r-1} \cup e}{\ul{e-r-1}\sm 0 \cup 2e};\]
  \item $G=\Spin_{2n}^+(q)$, $d$ is odd, $n = 2d+r$ with $0\leq r < d$ and
   \[\cS = \sy{\ul{d-r-1}\sm 0 \cup d \cup 2d}{\ul{d-r}};\]
  \item $G=\Spin_{2n}^+(q)$, $d=2e$ is even, $n=2e+r$ with $0\leq r \le e$ and
   \[\cS = \sy{\ul{e-r-1} \cup e}{\ul{e-r}\sm 0 \cup  2e};\]
  \item $G=\Spin_{2n}^-(q)$, $d$ is odd, $n = 2d+r$ with $0 < r \le d$ and
   \[\cS = \sy{\ul{d-r}\sm 0 \cup d \cup 2d}{\ul{d-r-1}};\]
  \item $G=\Spin_{2n}^-(q)$, $d=2e$ is even, $n=2e+r$ with $0 < r < e$ and
   \[\cS = \sy{\ul{e-r} \cup e}{\ul{e-r-1}\sm 0 \cup  2e};\]
  \item $G=\Sp_4(2)$, $d=2$ and $\cS =\sys{0\ 1\ 2}{-}$ or $\sys{0\ 2}{1}$
   (there is no regular semisimple $\ell$-singular element);
  \item $G=\Sp_6(2)$, $d=2$ and $\cS=\sys{0\ 1\ 3}{-}$;
  \item $G=\Sp_8(2)$, $d=4$ and $\cS=\sys{0\ 1}{4}$ or $\sys{1\ 4}{0}$; or
  \item $G=\Spin_8^-(2)$, $d=2,4$ and $\cS=\sys{1\ 3}{-}$.
 \end{enumerate}
 Conversely, all characters listed above take non-zero values on all regular
 semisimple $\ell$-singular classes. More precisely, in cases $G=\Spin_8^-(2)$,
 $d=2,4$, $\cS=\sys{1\ 3}{-}$ the character $\rho_\cS$ has value $2$ on classes
 of type $(-,(2,1,1))$. In all other cases the character values are $\pm 1$.
\end{thm}

\begin{proof}
(a) Our strategy is to use Corollary~\ref{cor:vanish} of the Murnaghan--Nakayama
formula. We consider maximal tori with $\ell$-singular regular
elements which are parameterised by pairs of partitions with very few
parts. For example, if for a symbol $\cS$ the character $\rho_\cS$ has
non-zero value on regular elements in a torus of type $((\la_1),(\mu_1))$
then the Murnaghan-Nakayama formula implies that we can remove a $\la_1$-hook
from $\cS$ and that the resulting symbol has a $\mu_1$-cohook which after
removing leaves a symbol of rank $0$.

(b) For groups of type $B_n$ and $C_n$ we need to consider symbols of rank
$n$ with odd defect. First note that there is only one symbol of rank $0$
with odd defect, namely the defect $1$ symbol
\[\sys{0}{-} = \sys{\ul{n}}{\ul{n-1}}\]
(apply the shift operation $n$ times). We will write
all symbols with $2n+1$ entries, because this way we can construct all
symbols of rank $n$ by adding hooks or cohooks starting from the rank $0$
symbol without the need of further shift operations. Also, we will write
all symbols of non-zero defect such that the larger number of entries
is in the first row.

The symbols of rank $t \leq n$ which have a $t$-hook are
$\sys{\ul{n}\sm k\cup k+t}{\ul{n-1}}$ for $n+1-t \leq k \leq n$ and their
dual symbols $\sys{\ul{n}}{\ul{n-1}\sm k\cup k+t}$ for $n-t\leq k\leq n-1$.

To find up to duality all symbols $\cS$ of rank $n$ which have an $(n-t)$-cohook
$c$ such that $\cS \sm c$ has a $t$-hook, it is sufficient to take the first
set of symbols of rank $t$ with a $t$-hook and to add an $(n-t)$-cohook in
all possible ways (filling the gap $k$, creating new entries in the first
or the second row or moving the $k+t$ entry). We get:
\[\begin{aligned}
  \sy{\ul{n}\cup l+n}{\ul{n-1}\sm l}&\quad (1\leq l\leq t),\\
  \sy{\ul{n}\sm k \cup k+t \cup l+n-t}{\ul{n-1}\sm l}&\quad (n-t< k \leq n,\
    t< l < n,\  k-l \neq n-2t),\\
  \sy{\ul{n-1}\cup l+n-t}{\ul{n}\sm k \sm l \cup k+t}&\quad
    (n-t< k \leq n,\ t \leq l \leq n, k\neq l),\\
  \sy{\ul{n-1}\cup k+n}{\ul{n}\sm k}&\quad (n-t< k \leq n).
\end{aligned}\]

(c) We now prove the theorem for the case $B_n$, $C_n$ for $n\geq 3$
and odd $d$. Write $n = ad+r$ with $0\leq r <d$. Note that for $q=2$ or $3$
we have $d > 1$ (since $\ell \neq 2$). In our argument we will use maximal tori
which by Lemma~\ref{lem:regelts} contain regular $\ell$-singular elements.

First assume $r=0$ and let $\cS$ be a symbol of rank $n$ such that $\rho_\cS$
is non-zero on all regular semisimple $\ell$-singular classes. Then as a
torus of type  $((n),-)$ contains $\ell$-singular elements, $\cS$ has
an $n$-hook, that is (up to duality) one of
\[ \cS = \sy{\ul{n}\sm k \cup k+n}{\ul{n-1}} \qquad (1\leq k\leq n).\]
Tori of type $((d),(n-d))$ also contain $\ell$-singular regular elements.
Therefore $\cS$ must have a $d$-hook $h$ and an $(n-d)$-cohook such that
$\cS\setminus h$ still has an $(n-d)$-cohook. This is only the case
for $k=n$ and $k=n-d$. For $k=n$ we obtain the symbol $\sys{\ul{n-1}\cup
2n}{\ul{n-1}}$ of the trivial character.  For $k = n-d$ we get the symbol
\[ \cS = \sy{\ul{n}\sm n-d \cup 2n-d}{\ul{n-1}}.\]
If $a>2$ then this symbol has no $(a-1)d = (n-d)$-hook, hence $\rho_\cS$ is
zero on regular elements in tori of type $((n-d),(d))$. For $a=2$ we get one
of the symbols shown in part~(1) of the statement (after applying the
inverse shift operation $n-d$ times).

Now let $r>0$ (and so $d\geq 3$). We consider tori of types
$((ad),(r)) = ((n-r),(r))$ and $((d,n-d),-)$. Here we check for every symbol
$\cS$ listed in~(b) (for $t=n-r$) if it has an $(n-r)$-hook and if
we can get the rank $0$ symbol by removing a $d$-hook and an $(n-d)$-hook
in both possible orders.
Note that the last condition implies that $\cS$ has defect one because
removing hooks does not change the defect. These conditions yield in all
cases relations between the parameters $k$ and $l$ and then determine $k$. The
only possible symbols which remain are those of the trivial character and
\[ \cS = \sy{\ul{n}}{\ul{n}\sm d \sm n-r \cup n+d-r}.\]
If $a>2$ this has no $(a-1)d$-hook and so $\rho_\cS$ is zero on regular
elements in tori of type $(((a-1)d),(d+r))$. For $a=2$ the dual of $\cS$ is
given in statement~(1), that is $\sys{\ul{n-1}\sm n-d+r\cup n+r\cup
n+d+r}{\ul{n-1}}$ and applying the inverse shift operation $n-d+r$ times.

(d) The argument for type $B_n$, $C_n$, $n \geq 7$, and even $d=2e$,
$n=ae+r$, is very
similar. Pairs of partitions corresponding to tori with $\ell$-singular regular
elements must now have an odd multiple of $e$ in the second part or an even
multiple of $e$ in the first part. Therefore we distinguish the cases
of even and odd $a$ and of $r=0$ and $r>0$. We just list the types of tori
which we consider in the given order such that the arguments from~(c) work
with slight modifications.

$a$ even, $r=0$: $((n),-)$ and $((n-e),(e))$, and for $a>2$ also
$((n-2e),(2e))$.

$a$ even, $r>0$: $((n-r),(r))$, $(-,(n-e,e))$ and $((n-e),(e))$, and for
$a>2$ also $(-,((a-1)e, e+r))$.

$a$ odd, $r=0$: $((n-e),(e))$, $(-,(n))$ and $(-,(n-e,e))$, and finally to
rule out a symbol of defect 3 also $((n-e,e),-)$ (if $a\leq 3$) or
$((n-3e), (3e))$ (if $a>3$).

$a$ odd, $r>0$: $((n-r-e),(e+r))$, $(-,(n-r,r))$ and $((n-e),(e))$.

(e) The cases of type $D_n$ and $\tw2D_n$ can also be discussed with very
similar arguments. The only symbol of rank $0$ and even defect has defect
$0$ and is $\sys{-}{-} = \sys{\ul{n-1}}{\ul{n-1}}$ and so all the symbols of
rank $n$ with an $n$-hook are
\[ \cS=\sy{\ul{n-1}\sm k \cup k+n}{\ul{n-1}} \quad (0\leq k \leq n).\]
Here it is sufficient to consider the cases with $k\geq (n-1)/2$ in the
further arguments because the remaining ones are the duals of these.

(f) Type $D_n$, $n\geq 7$. We list again which tori we consider in the
various cases.

$d$ odd, $r=0$:
$((n),-)$ and in case $d=1$ we use $((n-1, 1),-)$ and for $d\geq 3$
we use $((n-d),(d-1,1))$. For $a>2$ we also need $((n-2d),(2d-1,1))$ to rule
out one possibility.

$d$ odd, $r>0$: $((n-d,d),-)$ and $((n-r,r),-)$ if $r\neq 2$ or
$((n-2),(1,1))$ if $r=2$, $((d),(n-d-1,1))$.

$d=2e$ even, $n=ae$, $a$ even: $((n),-)$ and $(-,(n-1,1))$ if $e=1$
or $((n-e-1),(e,1))$ if $e>1$. For $a>2$ we also use $((n-2e),(2e-1,1))$.

$d=2e$ even, $n=ae+r$, $a$ even, $r>0$: $(-,(n-e,e))$, $((n-e-1),(e,1))$ and
$((n-r,r),-)$ if $r\neq 2$ or $((n-2),(1,1))$ if $r=2$. For $a>2$ we
also need $((2e),(n-2e-1,1))$.

$d=2e$ even, $n=ae$, $a$ odd: $(-,(n-e,e))$ and $((n-e,e),-)$ if $e \neq 2$
or $((n-2),(1,1))$ if $e=2$. Finally, use $((n-3e,3e),-)$ if $a > 3$. For
$a=3$ we get the exception listed in~(4) for the parameter $r=e$ (that is
$n = 2e+r = 3e$).

$d=2e$ even, $n=ae+r$, $a$ odd, $r>0$: $(-,(n-e,e))$, $(-,(n-r,r))$,
$((n-e-1),(e,1))$ and $((2e),(n-2e-1,1))$.

(g) Type $\tw2D_n$, $n\geq 7$. We list again which tori we consider in the
various cases.

$d$ odd, $n=ad$: Recall that in this case we have $a\geq
3$. We use $((n-d),(d))$ and $((n-2d),(2d))$. For $a > 3$ we also use
$((n-3d),(d))$ and for $a=3$ we get the exception listed in~(5) for $r=d$
(that is $n = 2d+r = 3r$).

$d$ odd, $n=ad+r$, $r>0$ (and so $d\geq 3$):
$((n-r),(r))$,  $((d),(n-d))$  and $((d,1),(n-d-1))$. And for $a>2$ also
$((2d),(n-2d))$. For $a=2$ this leads to the other exceptions listed in~(5).

$d=2e$ even, $n=ae$, $a$ even: Recall that in this case we
have $a\geq 4$. We use $((n-e),(e))$, $((n-2e),(2e))$ and one of
$((3),(n-3))$ (if $e=1$) or $((n-3e),(3e))$ (if $e>1$) or
$(-,(6,1,1))$ (if $n=8$).

$d=2e$ even, $n=ae+r$, $a$ even, $r>0$ (so $e>1$):
$((n-e),(e))$ and $((n-r),(r))$. For $a>2$ we also use $((2e),(n-2e))$.
For $a=2$ this leads to the exceptions listed in~(6).

$d=2e$ even, $n=ae$, $a$ odd: $(-,(n))$, $((n-e),(e))$ and
$((n-e-1,1),(e))$.

$d=2e$ even, $n=ae+r$, $a$ odd, $r>0$: $((n-e),(e))$, $((n-e-r),(e+r))$ and
$((2e),(n-2e))$. Finally also $((r),(n-r))$ if $r\neq 2$ or $(-,(n-2,1,1))$
for $r=2$.

(h) For small rank $n \leq 6$ we compute the number of regular elements in all
tori, and determine the values of unipotent characters on all
$\ell$-singular regular semisimple classes via the Murnaghan-Nakayama
formula in Theorem~\ref{thm:MN}. This yields the extra cases in statements~(7)
to~(10).

(i) The symbols given in statements~(1) to~(6) each contain a unique
$d$-hook and a unique $2d$-hook (or a unique $e$-cohook and unique
$2e$-hook, respectively). Removing these we always get the symbol of the
trivial representation of the appropriate rank (and the Steinberg
representation for the dual symbols). Therefore, the Murnaghan-Nakayama
formula shows that the values on all regular semisimple $\ell$-singular
elements are $\pm 1$.
\end{proof}

\section{Simple endotrivial modules}   \label{sec:endotriv}
In this section we prove Theorem~2 on the classification of simple endotrivial
modules for finite classical groups. Note that the case of linear and
unitary groups has already been treated in \cite{LM15}, so we only need to
consider symplectic groups and spin groups as in the previous section.

\subsection{On certain unipotent characters}   \label{subsec:unip}

We first rule out some candidate characters. Let $\Phi_d\in\ZZ[x]$ denote
the $d$th cyclotomic polynomial.

\begin{prop}   \label{prop:babbage}
 Let $\chi$ be the unipotent character of $G = \SO_{6d-1}(q)$ or $\Sp_{6d-2}(q)$
 (respectively of $\SO_{2(3d-1)}^+(q)$, $\SO_{6d}^-(q)$) parameterised by the
 symbol
 $$\binom{d\ 2d}{0}\quad \text{respectively }\ \binom{d\ 2d}{0\ 1},\
  \binom{d\ 2d}{},$$
 with $d\ge3$ odd, or the unipotent character of $\SO_{6e}^+(q)$
 parameterised by
 $$\binom{2e}{e}$$
 for $d=2e$ even. Then $\chi(1)\not\equiv\pm1\pmod{|G|_\ell}$ for
 any prime $2<\ell$ with $d_\ell(q)=d$.
\end{prop}

\begin{proof}
Since $\Phi_d(q)$ divides $|G|$ at least twice it is sufficient to show that
$\chi(1)\not\equiv\pm1\pmod{\Phi_d(q)_\ell^2}$.
We use a version of Babbage's congruence for quantum binomial coefficients
$$\left[hd-1\atop d-1\right]_x
   \equiv x^{(h-1)\binom{d}{2}}\pmod{\Phi_d(x)^2}\quad\text{in }\ZZ[x]$$
for any $d,h\ge2$ (proved in \cite[Lemma 3.7]{LM15}). If $\chi$ belongs to
$\binom{d\ 2d}{0}$, we can express its degree using a quantum binomial
coefficient evaluated at $q^2$, see~\cite[13.8]{Ca}:
$$\chi(1) = \frac{(q^2-1)\cdots(q^{6d-2}-1)(q^{2d}-q^d)(q^{2d}+1)(q^d+1)}
  {2(q^2-1)\cdots(q^{4d}-1)(q^2-1)\cdots(q^{2d}-1)}
  =\frac{1}{2}\left[3d-1\atop d-1\right]_{q^2}q^d(q^{2d}+1).$$
Substituting $x = q^2$ in the Babbage congruence we get a congruence over
the integers modulo $\Phi_d(q^2)^2$. Since $\Phi_d(x^2) =
\Phi_d(x)\Phi_d(-x)$ we also get a congruence modulo $\Phi_d(q)^2$:
\[\chi(1) \equiv \frac{1}{2} q^{2d(d-1)+d}(q^{2d}+1) \pmod{\Phi_d(q)^2}.\]
It suffices to show that $2\chi(1) \not\equiv \pm 2\pmod{\Phi_d(q)_\ell^2}$.
Since $\Phi_d(q) \mid (q^d-1)$ we get
\[\begin{array}{rcl}
  2\chi(1) &\equiv& q^{2d(d-1)+d} \cdot 2 q^d
   = 2 (q^{2d})^d \equiv 2(q^{2d}-(q^d-1)^2)^d = 2(2(q^d-1)+1)^d \\
  &\equiv& 2 (2d(q^d-1)+1) \pmod{\Phi_d(q)^2}.
\end{array}\]
Since $\ell \neq 2$, $d < \ell$  and $\ell \mid (q^d-1)$ we conclude that
$2\chi(1) + 2$ is not divisible by $\ell$ and that the highest power of
$\ell$ dividing $2\chi(1)-2$ is $(q^d-1)_\ell = \Phi_d(q)_\ell$.

Now we turn to the unipotent character with symbol $\binom{d\ 2d}{0\ 1}$.
Its degree differs from the previously considered $\chi(1)$ by the factor
$$\frac{(q^{2d}+q)(q^{d-1}+1)}{(q^{3d-1}+1)(q+1)}
  =\frac{q^{3d-1}+3q^d+q-1+(q^d-1)^2}{q^{3d-1}+3q^d+q-1+(q^d+2)(q^d-1)^2}.$$
We argue that both degrees are equal modulo $\Phi_d(q)_\ell^2$.

For this, note that for integers $a,b,c,m$ with $a b/c \in \ZZ$, $b \equiv c
\pmod{m}$ and $(c,m) = 1$ there is an integer $\tilde{c}$ with $c\tilde{c}
\equiv b \tilde{c}
\equiv 1 \pmod{m}$ and so $a b/c \equiv ab/c\cdot c\tilde{c} \equiv a b
\tilde{c} \equiv a \pmod{m}$. This applies here because $(q^{3d-1}+1)(q+1)$
is not divisible by $\ell$ (since $d\geq 3$ is the order of $q$ modulo
$\ell$ and $(q^{3d-1}+1)(q+1) (q^{3d-1}-1)(q-1) = (q^{6d-2}-1)(q^2-1)$ is
not divisible by $\ell$).

The same argument works for the degree of the unipotent character with symbol
$\binom{d\ 2d}{}$. Its degree differs from the previous $\chi(1)$ by the factor
\[\frac{2(q^{2d}-q^d+1)}{q^{2d}+1}
  = \frac{2q^d +2(q^d-1)^2}{2q^d+(q^d-1)^2}\]
and $q^{2d}+1$ is not divisible by $\ell$.

Finally we consider for even $d=2e$ the unipotent character $\psi$ with symbol
$\binom{2e}{e}$.   Its degree is
\[ \psi(1) = \left[3e-1\atop e-1\right]_{q^2} (q^{2e}+q^e+1)q^e.\]
Evaluating Babbage's congruence at $q^2$ yields a congruence modulo
$\Phi_e(q^2)^2$. But note that for odd $e$ we have $\Phi_e(x^2) = \Phi_e(x)
\Phi_e(-x)$ and $\Phi_e(-x) = \Phi_d(x)$, and for even $e$ we have
$\Phi_e(x^2) = \Phi_d(x)$. So, we also get congruences modulo $\Phi_d(q)^2$.
Using now  $\Phi_d(q) \mid (q^e+1)$ we find with similar arguments as
before that $\psi(1) \equiv -q^{2e^2} \pmod{\Phi_d(q)^2}$ and that
$\psi(1)-1$ is not divisible by $\ell$, and that $\psi(1)+1$ is divisible by
$\Phi_d(q)_\ell$ but no higher power of $\ell$.
\end{proof}

In order to deal with certain unipotent characters, we recall some
results of James and Mathas \cite{JM00} on reducibility of characters of Hecke
algebras. Let $\cH_n=\cH_n(q,Q)$ be the Iwahori--Hecke algebra of type $B_n$
with parameters $q$ and $Q$ over $\ZZ[q^{\pm1},Q^{\pm1}]$. Its irreducible
characters are in natural bijection with those of the Weyl group $W(B_n)$ via
the specialisation $q\mapsto1,Q\mapsto1$, and hence labelled by
pairs of partitions
$(\la,\mu)$ of $n$. We write $\chi_{\la,\mu}$ for the character of $\cH_n$
labelled by $(\la,\mu)$. We consider certain specialisations of $\cH_n$.

\begin{prop}   \label{prop:HeckeB}
 Let $\ell$ be a prime and $q$ a prime power. Set $d=d_\ell(q)$.
 \begin{enumerate}
  \item[\rm(a)] If $d$ is odd, $n=2d+r$ with $0\le r\le d-2$, and
   $(\la,\mu)=((d+r,r+1,1^{d-r-1}),-)$, then the character $\chi_{\la,\mu}$
   of $\cH_n(q,q)$ is reducible modulo~$\ell$.
  \item[\rm(b)] If $d=2e>2$ is even, $n=2e+r$ with $0\le r\le e-1$, and
   $(\la,\mu)=((r),(e+r+1,1^{e-r-1}))$, then the character $\chi_{\la,\mu}$
   of $\cH_n(q,q)$ is reducible modulo~$\ell$.
 \end{enumerate}
\end{prop}

\begin{proof}
Let $\zeta_d$ denote a primitive complex $d$th root of unity. The reduction
modulo~$\ell$ of $\cH_n(q,q)$ factors through the specialisation to
$\cH_n(\zeta_d,\zeta_d)$, thus a character which becomes reducible under the
latter specialisation will a fortiori be reducible modulo~$\ell$.
James and Mathas \cite[Thms.~4.7(i) and~4.10]{JM00} give a sufficient
criterion for $\chi_{\la,\mu}$ to become reducible under the specialisation to
$\cH_n(\zeta_d,\zeta_d)$, if $(\la,\mu)$ is a Kleshchev bipartition.
In case~(a), $\chi_{\la,\mu}$ is reducible modulo $\Phi_d$ by
\cite[Thm.~4.10(i)]{JM00} since the character of the Hecke algebra of $\fS_n$
labelled by $\la=(d+r,r+1,1^{d-r-1})$ is reducible modulo $\Phi_d$, as the
first column hook lengths are different. \par
In~(b), $(\la,\mu)$ is Kleshchev, and for $r>0$ the last box of the Young
diagram of $\la$ can be moved to the end of the first column of $\mu$ without
changing its residue, so $\chi_{\la,\mu}$ is reducible by
\cite[Thm.~4.10(ii)]{JM00}. For $r=0$ the character $\chi_{\la,\mu}$ is
reducible modulo $\Phi_{2e}$ by \cite[Thm.~4.10(i)]{JM00} since the
character of the Hecke algebra of $\fS_n$ labelled by the second part
$\mu=(e+1,1^{e-1})$ is.
\end{proof}

The next reducibility result for characters of specialisations of the Hecke
algebra of type $B_n$ will be used in the investigation of
unipotent characters of $\Spin_{2n}^\pm(q)$.

\begin{prop}   \label{prop:HeckeD}
 Let $\ell$ be a prime and $q$ a prime power. Set $d=d_\ell(q)$.
 \begin{enumerate}
  \item[\rm(a)] If $d$ is odd, $n=2d+r$ with $0\le r\le d-2$, and
   $(\la,\mu)=((d+r,r+1,1^{d-r-1}),-)$, then the character $\chi_{\la,\mu}$,
   of $\cH_n(q,1)$ is reducible modulo~$\ell$.
  \item[\rm(b)] If $d=2e>2$ is even, $n=2e+r$ with $0\le r\le e-1$, and
   $(\la,\mu)=((e+r,1^{e-r}),(r))$, then the character $\chi_{\la,\mu}$ of
   $\cH_n(q,1)$ is reducible modulo~$\ell$.
  \item[\rm(c)] If $d$ is odd, $n=2d+r$ with $1\le r\le d-1$, and
   $(\la,\mu)=(d+r-1,r,1^{d-r}),-$, then the character $\chi_{\la,\mu}$ of
   $\cH_{n-1}(q,q^2)$ is reducible modulo~$\ell$.
  \item[\rm(d)] If $d=2e>2$ is even, $n=2e+r$ with $1\le r\le e-1$, and
   $(\la,\mu)=((r-1),(e+r+1,1^{e-r-1}))$, then the character $\chi_{\la,\mu}$
   of $\cH_{n-1}(q,q^2)$ is reducible modulo~$\ell$ unless $(e,r)=(2,1)$.
 \end{enumerate}
\end{prop}

\begin{proof}
Parts~(a) and~(c) again follow from \cite[Thm.~4.10(i)]{JM00} as the first
column hook lengths of $\la$ are different. The characters
$\chi_{\la,\mu}$ as in~(b) and~(d) are reducible since the first, respectively
second part of $\la$ is a hook partition, of length~$d=2e$, hence belonging
to a reducible character of the Hecke algebra of type $\fS_{2e}$, unless this
hook is either $(2e)$ or $(1^{2e})$. The latter case only occurs in~(d) when
$(e,r)=(2,1)$.
\end{proof}

\subsection{Simple endotrivial modules in classical type groups}   \label{subsec:endoclass}

We will use the fact shown
in~\cite[Thm.~1.3]{LMS13} that any endotrivial module is liftable to a
characteristic~0 representation. Therefore, we can investigate endotrivial
modules by the ordinary character of their lift. These characters have the
following properties:
\begin{itemize}
 \item they are irreducible modulo $\ell$,
 \item their values on $\ell$-singular classes are of absolute value $1$,
 \item their degree is congruent to $\pm 1$ modulo $|G|_\ell$.
\end{itemize}
Theorem~2 is now a consequence of:

\begin{thm}   \label{thm:BC}
 Let $G=\Sp_{2n}(q)$ or $G=\Spin_{2n+1}(q)$ with $n\ge2$, $(n,q)\ne(2,2)$, or
 $G=\Spin_{2n}^\pm(q)$ with $n\ge4$, and $\ell$ a prime for which Sylow
 $\ell$-subgroups of $G$ are non-cyclic. Let $\rho$ be the character of a
 nontrivial simple endotrivial $kS$-module, for a central factor group $S$
 of $G$. Then we have $G=S=\Sp_8(2)\cong\OO_9(2)$, $\ell=5$, and
 $\rho=\rho_\cS$, with $\cS=\binom{0\ 1}{4}$, is the unipotent character of
 degree~$\rho_\cS(1)=51$.
\end{thm}

\begin{proof}
By \cite[Thms.~6.7 and~5.2]{LMS13} there do not exist non-trivial simple
endotrivial $kS$-modules if either $\ell=2$ or $\ell$ divides~$q$.
Now let $d=d_\ell(q)$. Since Sylow $\ell$-subgroups of $G$
are non-cyclic, the cyclotomic polynomial $\Phi_d(q)$ has to divide the group
order at least twice. By Theorem~\ref{thm:non-unip}, all non-unipotent
characters of $G$ vanish on some $\ell$-singular element of $G$, so cannot
come from simple endotrivial modules. The only candidates for endotrivial
unipotent characters are those appearing in the conclusion of
Theorem~\ref{thm:unip}. Of those, we can discard the Steinberg character
of $G$, because $G$ contains $\ell$-singular elements with non-trivial unipotent
part, and on those classes the value of the Steinberg character is zero.
\par
Now consider the other unipotent characters listed in Theorem~\ref{thm:unip}.
First assume that $G=\Sp_{2n}(q)$ or $G=\Spin_{2n+1}(q)$. Let $\rho=\rho_\cS$
with $\cS$ as in Theorem~\ref{thm:unip}(1). So $d$ is odd and $n=2d+r$,
$0\le r\le d-1$.  Then $\rho$
lies in the principal series and is parameterised by the character of the Weyl
group $W(B_n)$ labelled by the pair of partitions
$((d+r,r+1,1^{d-r-1}),-)\vdash 2d+r$. By Proposition~\ref{prop:HeckeB}(a), the
corresponding character of the Hecke algebra $\cH_n$ of type $B_n$ is
reducible modulo~$\ell$ for $r\ne d-1$. Since the decomposition matrix of the
Hecke algebra of $W(B_n)$ embeds into that of $G$ (see e.g.{}
\cite[Thm.~4.1.14]{GJ11}), this shows that $\rho_\cS$ is reducible
modulo~$\ell$, so does not give an example. The same argument applies to the
Alvis--Curtis dual characters, since Alvis--Curtis duality on the side of
Hecke algebras corresponds to tensoring with the sign character. In the case
$r=d-1$ we have $\cS=\binom{d\ 2d}{0}$ and the corresponding unipotent
character does not satisfy the necessary degree congruence by
Proposition~\ref{prop:babbage}.
\par
So now assume that $d=2e$ is even and write $n=2e+r$ with $0\le r\le e-1$.
Let $\rho=\rho_\cS$ as in Theorem~\ref{thm:unip}(2). Again, $\rho$ lies in
the principal series and is labelled by the pair of partitions
$((r+1),(e+r,1^{e-r-1}))\vdash 2e+r$. Hence it is reducible modulo $\ell$
by Proposition~\ref{prop:HeckeB}(b) if $e>1$. For $e=1$ we have $r=0$, $n=2$
and $\cS=\binom{0\ 1}{2}$. Here
$\rho_\cS(1)=q(q^2+1)/2\equiv 2q+1\pmod{(q+1)^2}$ is not congruent to~$\pm1$
modulo the $\ell$-part of the group order. \par

For $G=\Sp_{2n}(q)$ or $G=\Spin_{2n+1}(q)$ it remains to consider the unipotent
characters listed in Theorem~\ref{thm:unip}(8) and~(9). (Note that the group
$\Sp_4(2)$ was excluded.) In case~(8)
we have $q=d=2$, so $\ell=3$. The listed symbol parameterises a unipotent
character of degree~7, its Alvis-Curtis dual has degree~56, both incongruent
$\pm1\pmod{3^3}$. Finally in case~(9), $\ell|(q^2+1)=5$. The second symbol
belongs to a character of degree~119, its Alvis--Curtis dual has degree~30464,
both incongruent to $\pm1\pmod{5^2}$, so this gives no example.
But the 51-dimensional simple unipotent module for $\Sp_8(2)$ labelled by
$\binom{0\ 1}{4}$ is irreducible and endotrivial in characteristic $\ell=5$.
Indeed, the restriction of the corresponding irreducible character to the
subgroup $\OO_8^+(2)$ contains the trivial character once and all other
constituents are of 5-defect zero, hence projective modulo~5. Its Alvis--Curtis
dual has degree~$13056\equiv6\pmod{5^2}$, so cannot lead to an example.
\par\medskip
Next consider the case that $G=\Spin_{2n}^+(q)$. Let $\rho=\rho_\cS$ be
unipotent with $\cS$ as in
Theorem~\ref{thm:unip}(3). So $d$ is odd and $n=2d+r$ with $0\le r\le d-1$.
Then $\rho$ lies in the principal series and is parameterised by the character
of the Weyl group $W(D_n)$ labelled by the pair of partitions
$(\la,\mu)=((d+r,r+1,1^{d-r-1}),-)\vdash 2d+r$. By
Proposition~\ref{prop:HeckeD}(a), the corresponding character $\chi_{\la,\mu}$
of the Hecke algebra $\cH_n$ of type $B_n$ is reducible modulo~$\ell$ for
$r\ne d-1$. Since $\la\ne\mu$, $\chi_{\la,\mu}$ restricts irreducibly to the
Hecke algebra $\cH_n'$ of type $D_n$ (of index~2) and thus is also reducible as
a character of $\cH_n'$ modulo $\ell$. Using that the decomposition matrix of
the Hecke algebra of $W(D_n)$ embeds into that of $G$, this shows that
$\rho_\cS$ is reducible modulo~$\ell$, so does not give an example. The same
argument applies to the Alvis--Curtis dual characters. In the remaining case
$r=d-1$, Proposition~\ref{prop:babbage} shows that the corresponding unipotent
character does not satisfy the degree congruence.
\par
Now assume that $\rho_\cS$ is as in Theorem~\ref{thm:unip}(4), so $d=2e$ is
even and $n=2e+r$ with $0\le r\le e$. Again, $\rho$ lies in the principal
series and is labelled by the pair of partitions
$((e+r,1^{e-r}),(r))\vdash 2e+r$. By Proposition~\ref{prop:HeckeD}(b) the
corresponding character of the Hecke algebra of type $B_n$ is reducible
modulo~$\ell$ if $r<e$ and $e>1$, and so as before the same holds for the
character of the Hecke subalgebra of type $D_n$, whence $\rho_\cS$ is reducible
modulo~$\ell$ for $e>1$. Note that $e=1$ is excluded since then $n=2e=2$. The
case when $r=e$ does not satisfy the degree congruence by
Proposition~\ref{prop:babbage}.
\par
Similarly, in the case that $G=\Spin_{2n}^-(q)$ the characters of the principal
series Hecke algebra $\cH_{n-1}(q,q^2)$ labelled by $(\la,\mu)$ as in
Theorem~\ref{thm:unip}(5) and~(6) are reducible modulo~$\ell$ by
Proposition~\ref{prop:HeckeD}(c) and~(d), so the same holds for the unipotent
characters parameterised by these pairs, unless $r=d$. The latter character
again does not belong to an endotrivial module by
Proposition~\ref{prop:babbage}. In the exceptional case $(e,r)=(2,1)$ of
Proposition~\ref{prop:HeckeD}(d) we have $n=2e+r=5$
and $\rho_\cS(1)+1\equiv 4(q^2+1)\pmod{(q^2+1)^2}$, hence the necessary degree
congruence is not satisfied for $\rho_\cS$.

Finally, the character listed in Theorem~\ref{thm:unip}(10) takes value~2 on
some $\ell$-singular classes, by the final statement in that result, so again
cannot be endotrivial.
\end{proof}

The endotrivial characters of $\Sp_4(2)\cong\fS_6$
and its covers have been described in \cite[Thm.~4.9]{LMS13}.

As pointed out before the simple endotrivial modules of linear and unitary
groups have already been determined in \cite[Thm.~3.10 and~4.5]{LM15}.

\section{Zeroes of characters of quasi-simple groups }   \label{sec:rank3}

Together with our previous work \cite{LM15,LMS13} with C.~Lassueur and
E.~Schulte our results allow us to guarantee the existence of zeroes of
characters of quasi-simple groups on $\ell$-singular elements once the
$\ell$-rank is at least~3, as claimed in Theorem~\ref{thmA:zeroes} which we
restate:

\begin{thm}   \label{thm:zeroes}
 Let $\ell>2$ be a prime and $G$ a finite quasi-simple group of $\ell$-rank at
 least~3. Then for any non-trivial character $\chi\in\Irr(G)$ there exists an
 $\ell$-singular element $g\in G$ with $\chi(g)=0$, unless one of:
 \begin{enumerate}
  \item[\rm(1)] $G$ is of Lie type in characteristic~$\ell$;
  \item[\rm(2)] $\ell=5$, $G=\PSL_5(q)$ with $5||(q-1)$ and $\chi$ is unipotent
   of degree $\chi(1)=q^2\Phi_5$;
  \item[\rm(3)] $\ell=5$, $G=\PSU_5(q)$ with $5||(q+1)$ and $\chi$ is unipotent
   of degree $\chi(1)=q^2\Phi_{10}$; or
  \item[\rm(4)] $\ell=5$, $G=Ly$ and $\chi(1)\in\{48174,11834746\}$; or
  \item[\rm(5)] $G=E_8(q)$ with $q$ odd, $d_\ell(q) = 4$ and $\chi$ is the 
   semisimple character with label $((0),(8))$ in the unique 
   Lusztig-series of type $D_8$.
 \end{enumerate}
\end{thm}

\begin{proof}
For $G$ a covering group of a sporadic simple group, an easy check of the
known character tables shows that only case~(4) arises. For alternating groups
and their covering groups the claim is contained in the proof of
\cite[Prop.~4.2 and Thm.~4.5]{LMS13}. If $G$ is a covering group of a
special linear or unitary group then \cite[Thm.~1.3]{LM15} shows that only the
cases in~(1), (2) and~(3) can arise.  \par
For the remaining groups of Lie type first note that the Steinberg character
vanishes on the product of any $\ell$-element with a unipotent element
in its centraliser, so can be discarded from our discussion. We now first
consider groups $G$ of exceptional Lie type. Since the $\ell$-rank of $G$ is
at least~3, it must be of type $F_4$, $E_6,\tw2E_6$, $E_7$ or $E_8$. For
these it is shown in the proof of \cite[Thm.~6.11]{LMS13} that their irreducible
characters $\chi$ vanish on some $\ell$-singular element unless
either $\chi$ is unipotent or $\ell|(q^2+1)$ in $G=E_8(q)$ and $\chi$ lies in
the Lusztig series of an isolated element with centraliser $D_8(q)$. We deal
with these cases in turn.  \par
We begin with $G=E_8(q)$. Since $G$ has $\ell$-rank at least~3 we must have
$d=d_\ell(q)\in\{1,2,3,4,6\}$. First assume that $d=1$.
There exist maximal tori of orders divisible by $\Phi_1\Phi_7$, $\Phi_1\Phi_9$
and $\Phi_1\Phi_{14}$, hence $G$ contains semisimple $\ell$-singular elements
of order divisible by a Zsigmondy prime divisor $r$ of $\Phi_7,\Phi_9$ and
$\Phi_{14}$. All non-trivial unipotent characters apart from
$\phi_{8,1},\phi_{8,91}$ and the Steinberg character are of defect
zero for one of these three primes, hence vanish on any $r$-singular element.
For the remaining two characters, direct calculation with the character formula
shows that they vanish on elements of order $\Phi_1\Phi_4$. Very similar
considerations apply to the remaining cases $d\in\{2,3,4,6\}$. If $d=4$ and
$\chi$ lies in the Lusztig series of an isolated element with centraliser
$D_8(q)$, the same argument using Zsigmondy prime divisors of $\Phi_5$,
$\Phi_8$, $\Phi_{10}$ and $\Phi_{12}$ only leaves the four characters in that
series, those with label $\sys{0}{8}$ and $\sys{0 1}{1 8}$ and their duals.
Three of them are zero on suitable products of an $\ell$ element with a
unipotent element.
\par
The arguments for the unipotent characters of the remaining exceptional groups
are similar and easier. \par
Finally, let us suppose that $G$ is of classical Lie type. By
Theorem~\ref{thm:non-unip} the claim holds if $\chi$ is not
unipotent. For unipotent characters $\chi$, the claim follows from
Theorem~\ref{thm:unip}, since all exceptions listed there have $\ell$-rank
two, including those in cases~(4) and~(5) with $n=3e$ respectively $n=3d$.
\end{proof}

\begin{rem}
(a) Let $p$ be a prime and $f\ge1$. Then $G=\SL_2(p^f)$ has $p$-rank~$f$, but
its irreducible characters of degree $q\pm1$ do not vanish on any $p$-singular
elements of $G$. Similarly, there exist such characters for $\SL_3(p^f)$,
$\SU_3(p^f)$ and $^2G_3(3^{2f+1})$. This shows that Case~(1) in
Theorem~\ref{thm:zeroes} is a true exception. Also, $Ly$ has characters of
degrees 48174 and 11834746 which do not vanish on 5-singular elements. (It is
well-known that the Lyons group behaves like a characteristic~5 group in many
respects, see e.g.~\cite{PR}.)
It can be checked that the unipotent characters of $\PSL_5(q)$ and $\PSU_5(q)$
listed in Cases~(2) and~(3) are also true exceptions. \par
(b) We expect the exceptions in Case~(1) to only occur for groups of Lie
type in characteristic~$\ell$ of small Lie rank.  \par
(c) The results of \cite{LMS13}, \cite{LM15} and of the present paper show
that even for $\ell$-rank~2 there exist only relatively few irreducible
characters of quasi-simple groups not vanishing on some $\ell$-singular
element. Nevertheless we refrain from attempting to give an explicit list in
that case.
\end{rem}

\section{An application to small 1-PIMs }   \label{sec:PIMs}

Our results so far may also be used in order to investigate the first Cartan
invariant of finite simple classical groups. In their recent paper \cite{KKS},
Koshitani, K\"ulshammer and Sambale raised the question to understand simple
groups with non-cyclic Sylow $\ell$-subgroup for which the first Cartan
invariant $c_{11}$ equals~2.

\begin{thm}   \label{thm:cartan}
 Let $G$ be a finite simple group of classical Lie type $B_n,C_n,D_n$ or
 $\tw2D_n$. Let $\ell>2$ be a prime for which Sylow $\ell$-subgroups of $G$
 are non-cyclic. Then the $\ell$-modular projective cover of the trivial
 character of $G$ has at least three ordinary constituents. In particular
 $c_{11}\ge3$.
\end{thm}

\begin{proof}
For $\ell$ the defining characteristic of $G$, this was already shown in
\cite{KKS}. Now assume that $\ell \neq 2$ is different from the defining
characteristic of $G$.
Assume that the character of the $\ell$-modular projective cover of the
trivial $G$ module has the form $1_G+\chi$, with $\chi\in\Irr(G)$. Since
projective characters vanish on all $\ell$-singular elements, this implies
that $\chi(g)=-1$ on all elements $g\in G$ of order divisible by $\ell$.

We now use our results on zeroes of characters. By Theorem~\ref{thm:non-unip}
any non-unipotent character of $G$ vanishes on some $\ell$-singular element.
(Again note that $\Sp_4(2)\cong\fS_6$ is not simple.)
So $\chi$ must be unipotent. In this case, Theorem~\ref{thm:unip}
again shows that most unipotent characters also vanish on some
$\ell$-singular element. It remains to consider the listed exceptions.

The trivial character never has value $-1$. To rule out the Steinberg
character we use that there is always a non-regular element $s \in G$
of order $\ell$, hence there is a unipotent element $u \neq 1$ in the
centraliser of $s$. The Steinberg character vanishes on the non-semisimple
element $su$. For the small cases listed in~\ref{thm:unip}(7) to~(10)
and their Alvis--Curtis duals we just check that there is always an
$\ell$-singular regular  semisimple element on which the character has
value $\neq -1$.

It remains to consider the unipotent characters corresponding to symbols
in Theorem~\ref{thm:unip}(1) to~(6) and their
duals. We have mentioned that the values of these characters on
$\ell$-singular regular semisimple classes are $\pm 1$. We give for each
case a pair of parameters of maximal tori on which character values have
opposite signs. This is easy to compute with the Murnaghan--Nakayama
formula~\ref{thm:MN} using the definitions of the signs occurring there.
It turns out that in each case the same pair of tori works for the symbol
and its dual. The dual symbols and the pairs of tori are as follows:
\begin{itemize}
\item[(1)] $\cS^\vee = \sys{\ul{d+r}}{\ul{d+r}\sm 0 \sm d \cup 2d}$, tori
$((d),(d+r))$ and $((2d),(r))$,
\item[(2)] $\cS^\vee = \sys{\ul{e+r}\sm 0 \cup 2e}{\ul{e+r}\sm e}$, tori
$((2e),(r))$ and $((e+r),(e))$,
\item[(3)] $\cS^\vee = \sys{\ul{d+r-1}}{\ul{d+r}\sm 0 \sm d \cup 2d}$, tori
$((d,d+r),-)$ (or $((d,d-2),(1,1)$ if $r=0$) and
$((2d,r),-)$ (or $((2d),(1,1)$ if $r=2$),
\item[(4)] $\cS^\vee = \sys{\ul{e+r-1} \sm 0 \cup 2e}{\ul{e+r}\sm e}$, tori
$((2e),-)$ and $((e-1),(e,1))$ if $r=0$ (use $((1),(3,2))$ for $d=n=6$)
and in case $r>0$  tori
$(-,(e+r,e))$ and $((2e,r),-)$ (if $r\neq 2$) or $((2e),(1,1))$ (if $r = 2$),
\item[(5)] $\cS^\vee = \sys{\ul{d+r}}{\ul{d+r-1}\sm 0 \sm d \cup 2d}$, tori
$((d),(d+r))$ and $((2d),(r))$,
\item[(6)] $\cS^\vee = \sys{\ul{e+r}\sm 0 \cup 2e}{\ul{e+r-1} \sm e}$, tori
$((e+r),(e))$ and $((2e),(r))$.
\end{itemize}
The only cases which are not covered by this argument because one of the
given tori has no regular element are: $\Sp_4(2)$ with $d=2$ ($\ell=3$) and
$\Sp_8(2)$ with $d=4$ ($\ell=5$). The first group is not simple, and for the
second the claim can be checked from its character
table.
\end{proof}


\end{document}


\title[Corrigendum: Murnaghan--Nakayama rule]
      {Corrigendum: A Murnaghan--Nakayama rule\\ for values of unipotent
       characters in\\ classical groups}

\date{\today}

\author{Frank L\"ubeck}
\address{Lehrstuhl D f\"ur Mathematik, RWTH Aachen, Pontdriesch 14/16,
  52062 Aachen, Germany.}
\email{Frank.Luebeck@math.rwth-aachen.de}

\author{Gunter Malle }
\address{FB Mathematik, TU Kaiserslautern, Postfach 3049,
         67653 Kaisers\-lautern, Germany.}
\email{malle@mathematik.uni-kl.de}




\begin{abstract}
We settle a missing case in the proof of one of the main applications of
our results in \cite{LM16}.
\end{abstract}

\maketitle

In \cite{LM16} we derived a Murnaghan--Nakayama formula for the values of
unipotent characters of finite classical groups. As one application, we showed
that the first Cartan invariant $c_{11}^{(\ell)}$ of a finite simple classical
group $G$ in non-defining characteristic~$\ell$ is never equal to~2 if Sylow
$\ell$-subgroups of $G$ are
non-cyclic. But contrary to what was claimed in the first line of the proof
of \cite[Thm.~7.1]{LM16}, the same assertion for the case of defining
characteristic does not follow from the work of Koshitani, K\"ulshammer and
Sambale \cite{KKS}. We thank Shigeo Koshitani and J\"urgen M\"uller for
pointing this out to us.

Here, we give a short argument for the missing case. The argument also
covers the
analogous statement in \cite[Prop.~6.3]{LM15} pertaining to simple groups of
exceptional Lie type. The crucial observation is that Sylow $p$-subgroups in
groups of Lie type of characteristic~$p$ are large. Except for a few small
rank cases we do not need any knowledge about the ordinary character table of
the group $G$. 

Note that if $c_{11}^{(p)}=2$ then there exists $\chi\in\Irr(G)$ such that
$1_G+\chi$ is the character of a projective module; in particular $\chi(g)=-1$
for all $p$-singular elements $g\in G$.

\begin{lem}   \label{lem:main}
 Let $G$ be a finite group and $p$ a prime. Assume that $H\le G$ is a
 subgroup such that $|G:H|$ is divisible by $p$, but $|G:H|<|G|_p$. Then
 $c_{11}^{(p)}>2$.
\end{lem}

\begin{proof}
The permutation character $\pi$ of $G$ on the cosets of $H$ has exactly one
trivial constituent $1_G$. On the other hand, as $|G:H|$ is divisible by $p$,
the modular reduction of the permutation representation has (at least) two
trivial modular constituents, one in the socle, one in the head. Thus the
$p$-modular reduction of one of the non-trivial ordinary constituents $\chi$
of $\pi$ must contain the trivial $p$-modular character of $G$. By Brauer
reciprocity,
$\chi$ then occurs in the projective cover $\Psi$ of the trivial character. If
we had $c_{11}^{(p)}=2$ then necessarily $\Psi=1_G+\chi$. Since the degree of
any projective character is divisible by $|G|_p$, this is a contradiction to
$$|G|_p\le\Psi(1)=1_G(1)+\chi(1)\le\pi(1)=|G:H|.$$
\end{proof}

\begin{prop}   \label{main}
 Let $G$ be a finite simple group of Lie type in characteristic~$p$. Then
 $c_{11}^{(p)}>2$ unless $G=\PSL_2(p)$.
\end{prop}

\begin{proof}
If the Lie-rank of $G$ is not very small, then we can find a subgroup of $G$
which fulfills the condition of Lemma~\ref{lem:main}. In Table~\ref{tab} we
give a subgroup $\hat H$ (see e.g.~\cite[\S13]{MT} for the existence) of a
suitable central extension $\hat G$ of $G$.  From the order formulae for 
these groups (see for example~\cite[Tab.~24.1]{MT})
it is easy to check that $\hat H \leq \hat G$ fulfills the condition of
Lemma~\ref{lem:main}. Since the center of $\hat G$ is of order prime to $p$
the same holds for $H \leq G$ where $H$ is the image of $\hat H$ in $G$.

It remains to consider the groups of type $A_1$, $A_2$, $\tw2A_2$,
$\tw3D_4$, $\tw2B_2$, $\tw2G_2$ and $\tw2F_4$.  For all of these the
(generic) character tables are known and available in \Chevie~\cite{Chv}. 
Here we look for all irreducible characters $\chi$ which have value $-1$ on all
$p$-singular (non-semisimple) elements, so that $\chi+1_G$ is a candidate
for the $1$-PIM of $G$.
It is easy to see that apart from two exceptions all characters either
vanish on some $p$-singular element or that a value on a non-trivial
unipotent class is positive or that the values on two non-trivial 
unipotent classes are different.

The two exceptions occur in type $A_1$, so $G = \PSL_2(q)$, which has a
family of characters of degree $q-1$ with value $-1$ on all $p$-singular
elements,
and in type $A_2$, so $G = \PSL_3(q)$, which has a family of characters
of degree $q^3-1$ which also take value $-1$ on all $p$-singular elements.
By \cite[Thm.~5.2]{MW08} the $1$-PIM of $G=\PSL_2(q)$ has degree $(2^f-1)q$,
where $q=p^f$. Since the maximal character degree for $G$ is $q+1$, this
shows that $c_{11}^{(p)}>2$ unless $f=1$, that is, unless $q=p$ is prime, in
which case the Sylow $p$-subgroups of $G$ are cyclic (and so indeed
$c_{11}^{(p)}=2$). 
Similarly, by~\cite[Thm.~5.3]{MW08} the $1$-PIM of $G=\PSL_3(q)$, $q=p^f$,
has degree $(12^f-6^f+1)q^3$, so none of the candidate $1_G+\chi$ mentioned
above can be the $1$-PIM.
\end{proof}

\begin{table}[htbp]
 \caption{Large subgroups}   \label{tab}
\[\begin{array}{|llc||cc|}
\hline
     \hat G&& \hat H&  \hat G& \hat H\cr
\hline
 \SL_4^\pm(q)&& \Sp_4(q)&                   G_2(q)& \tw2A_2(q)\\
 \SL_n^\pm(q)& (n\ge5)& \GL_{n-1}^\pm(q)&   F_4(q)& B_4(q)\\
 \Spin_{2n+1}(q)& (n\ge2)& \Spin_{2n}^-(q)&    E^\pm_6(q)& F_4(q)\\
 \Sp_{2n}(q)& (n\ge3)& \Sp_{2n-2}(q)\Sp_2(q)&      E_7(q)& E_6(q)\\
 \Spin_{2n}^\pm(q)& (n\ge4)& \Spin_{2n-1}(q)&  E_8(q)& E_7(q)\\
\hline
\end{array}\]
\end{table}

Observe that the case $p=2$ was already dealt with in \cite[Thm.~4.5]{KKS}.
